\documentclass[11pt]{article}
\usepackage{}

\usepackage[left=1in, right=1in, top=1in, bottom=1in]{geometry}
\usepackage{amsfonts}
\usepackage{amsthm}
\usepackage{amssymb}
\usepackage{amsmath,amsfonts}
\usepackage{mathrsfs}
\usepackage{cases}
\usepackage{latexsym,bm}
\usepackage{indentfirst}
\usepackage{color}
\usepackage{ifpdf}
\usepackage{graphicx}
\usepackage{psfrag}
\usepackage{dsfont}
\usepackage[
pdfauthor={CSY},
pdftitle={Strong generalized Halin graphs},
pdfstartview=XYZ,
bookmarks=true,
colorlinks=true,
linkcolor=blue,
urlcolor=blue,
citecolor=blue,
bookmarks=true,
linktocpage=true,
hyperindex=true
]{hyperref}

\newcommand{\xiaowuhao}{\fontsize{9pt}{\baselineskip}\selectfont}

\newtheorem{THM}{\textbf{Theorem}}[section]

\newtheorem{LEM}{\textbf{Lemma}}[section]

\newtheorem{COR}{\textbf{ Corollary}}[section]
\newtheorem{PRO}{\textbf{Proposition}}
\newcommand{\pf}{\textbf{Proof}.\quad}

\newtheorem{FAC}{\textbf{Fact}}

\newtheorem{CLA}{\textbf{Claim}}[section]
\newcommand{\qqed}{\hfill $\blacksquare$\vspace{1mm}}

\newcommand{\ve}{\varepsilon }

\newcommand{\scr}{\scriptscriptstyle}

\begin{document}

\title{Minimum degree condition for spanning generalized Halin graphs }
\author{ Guantao Chen, Songling Shan, and Ping Yang
\\{\xiaowuhao  Georgia State University, Atlanta, GA\,30303, USA}}
\date{}
\maketitle
\emph{\textbf{Abstract}.}
A spanning tree with no vertices of degree 2 is called a Homeomorphically
irreducible spanning tree\,(HIST).  Based on a HIST embedded in the plane,
a Halin graph is formed by connecting the leaves of the
tree into a cycle following the cyclic order determined by the embedding.
Both of the determination problems of whether a graph contains a HIST or
whether a graph contains a
spanning Halin graph are shown to be NP-complete.
It was
conjectured by Albertson, Berman, Hutchinson, and Thomassen in
1990 that a {\it every surface triangulation of at least four vertices contains a HIST}\,(confirmed).
And it was conjectured by
Lov\'asz and Plummer that {\it every 4-connected plane
triangulation contains a spanning Halin graph}\,(disproved).
Balancing the above two facts,
in this paper, we consider generalized Halin graphs, a family of graph structures which are
``stronger'' than HISTs but ``weaker'' than  Halin graphs in the sense of their construction constraints.
To be exact,
 a generalized Halin graph is formed from a HIST by connecting its leaves into a cycle.
Since a generalized Halin graph needs not to be planar, we investigate the
minimum degree condition for a graph to contain it as a spanning subgraph.
We show that there exists a positive integer $n_0$ such that any 3-connected graph with $n\ge n_0$ vertices and
minimum degree at least $(2n+3)/5$ contains a spanning generalized Halin graph.
As an application, the result implies that under the same condition,
the graph $G$ contains a wheel-minor of order at least $n/2$.
The minimum degree condition in the result is best possible.


\vspace{2mm}

\section{Introduction}

A tree with no vertex of degree 2 is called a {\it homeomorphically irreducible tree\,(HIT)},
and a spanning tree with no vertex of degree 2 is a {\it homeomorphically irreducible spanning tree\,(HIST)}.
A {\it Halin graph}, constructed by Halin in 1971~\cite{Halin-halin-graph},
is a  graph formed from a plane embedding of a HIST by connecting its leaves
into a cycle following the cyclic order determined by the embedding.
In 1990,
Albertson, Berman, Hutchinson, and Thomassen~\cite{ABHT}
showed that {\it it is NP-complete to determine whether a graph contains a HIST}.
However, for special graph classes such as triangulations
of surfaces, they
conjectured that {\it every triangulation of
a surface with at least 4 vertices contains a HIST}.
The conjecture was confirmed in~\cite{HIST1}.
It was shown by Horton, Parker, and Borie~\cite{Halin-NP}
that {\it it is NP-complete to determine whether a graph contains a (spanning) Halin graph}.
Again, restricted to triangulations,
Lov\'asz and Plummer~\cite{LP-Halingraph-Con} conjectured that {\it every 4-connected plane
triangulation contains a spanning Halin graph}.
But the conjecture was disproved recently~\cite{Disprove-LP-Con}.
Since a Halin graph possesses many
hamiltonian properties\,(e.g., see~\cite{Bondy-pancyclic,Cornuejols-halin,BondyLov-halin-panciclicity}),
it seems that a graph has to have very ``good properties'' in order to
contain a Halin graph as a spanning subgraph. For this reason,
by relaxing on the planarity requirement, we
define a {\it generalized Halin graph} as
a graph formed from a HIST by connecting its leaves
into a cycle, and we study sufficient conditions
for implying the containment of
a spanning generalized Halin graph in a given graph.

Compared to Halin graphs,
generalized Halin graphs are less studied.
Kaiser et al.
in~\cite{Halin-prism-hamiltonian}
showed that a generalized Halin graph is
{\it prism Hamiltonian}; that is, the Cartesian product of
a generalized Halin graph and $K_2$ is hamiltonian.
Since a tree with no degree 2 vertices has
more leaves than the non-leaves, a generalized Halin
graph contains a cycle of length at least half of
its order. Also, one can notice
that by contracting the non-leaves of the underlying tree of a generalized
Halin graph into a singe vertex, a wheel graph is resulted
with the contracted vertex as the hub,
where a minor of a graph is
 obtained from the graph by deleting edges/contracting edges, or deleting vertices.
Therefore, a generalized Halin graph contains a wheel-minor of
order at least half of its order.
The investigation on the
properties of generalized Halin graphs is not of
the interest of this paper. Instead, in this paper,
we show the following two results.

\begin{THM}\label{Main-Result1}
It is NP-complete to determine whether a graph contains a spanning generalized Halin graph. 	
\end{THM}

\begin{THM}\label{Main-Result2}
There exists a positive integer $n_0$ such that every 3-connected graph with $n\ge n_0$ vertices and
minimum degree at least $(2n+3)/5$ contains a spanning generalized Halin graph. The result is best possible
in the sense of the connectivity and minimum degree constraints. 	
\end{THM}

Since a generalized Halin graph of order $n$ contains a wheel-minor
of order at least $n/2$, we get the following corollary.

\begin{COR}
There exists a positive integer $n_0$ such that every 3-connected graph with $n\ge n_0$ vertices and
minimum degree at least $(2n+3)/5$ contains a wheel-minor of order at least $n/2$. 	
\end{COR}

For notational convenience,
for a graph $T$, we denote by $L(T)$ the set of degree 1 vertices of $T$ and $S(T)=V(T)-L(T)$.
Also we abbreviate {\it spanning generalized Halin graph} as {\it SGHG} in what follows,
and denote a generalized Halin graph as $H=T\cup C$, where $T$ is the underlying HIST
of $H$ and $C$ is the cycle spanning  on $L(T)$.
The remaining of the paper is organized as follows. In Section 2, we prove Theorem~\ref{Main-Result1}
and show the sharpness of Theorem~\ref{Main-Result2}. In Section 3, we introduce some
notations and lemmas, which are used in the proof of Theorem~\ref{Main-Result2}.
We then proof Theorem~\ref{Main-Result2} in Section 4.

%

\section{Proof of Theorem~\ref{Main-Result1} and the sharpness of Theorem~\ref{Main-Result2}}

\proof[Proof of Theorem~\ref{Main-Result1}]
It was shown by Albertson et al.~\cite{ABHT} that
it is NP-complete to decide whether a graph contains a HIST,
and by the definition, a generalized Halin graph contains a
HIST. Hence,
we see that the problem of deciding whether an
arbitrary graph contains an SGHG is in NP.
To show the problem is NP-complete we assume
the existence of a polynomial algorithm to test for an SGHG and
use it to create a polynomial algorithm to test for
a hamiltonian path between two vertices in an arbitrary graph. The
decision problem for such hamiltonian paths is a classic NP-complete
problem~\cite{NP-completeness}.

Let $G$ be a graph and  $x, y\in V(G)$. We want to determine
whether there exists a hamiltonian path connecting $x$ and $y$.
We first construct a new graph $G'$ and show that $G$
contains a hamiltonian path between $x$ and $y$ if and only if
$G'$ contains a HIST\,(the proof of this part is the same as the proof of Albertson et al. in~\cite{ABHT}).
Then based on $G'$, we construct a graph $G''$ and show that $G'$ contains
a HIST if and only if $G''$ contains an SGHG.

Let $\{z_1,z_2,\cdots, z_t\}=V(G)-\{x,y\}$. Then $G'$ is
formed by adding new vertices $\{z_1', z_2',\cdots, z_t'\}$ and
new edges $\{z_iz_i'\,:\, 1\le i\le t\}$. It is clear
that if $P$ is a hamiltonian path between $x$ and $y$,
then $P\cup \{z_iz_i'\,:\, 1\le i\le t\}$ is a HIST of $G'$.
Conversely, let $T$ be a HIST of $G'$. Since $1\le d_T(z_i')\le d_{G'}(z_i')=1$,
we get $d_T(z_i')=1$ for each $i$. Since $N_{G'}(z_i')=\{z_i\}$ and
$T$ is a HIST, we have $d_T(z_i)\ge 3$.  Hence $T-\{z'_1,z'_2,\cdots, z'_t\}$
is a tree with leaves possibly in $\{x,y\}$.
Since each tree has at least 2 leaves and a tree with exactly two leaves
is a path, we conclude that $T-\{z'_1,z'_2,\cdots, z'_t\}$ is a path between
$x$ and $y$.

Then based on $G'$, we construct a graph $G''$. First we add new vertices
$\{z'_{i1}, z_{i2}', z_{i3}'\,:\, 1\le i\le t\}$. Then we add edges
$\{z_i'z'_{i1}, z_i'z_{i2}', z_i'z_{i3}', z'_{i1}z_{i2}', z_{i2}'z_{i3}'\,:\, 1\le i\le t\}$.
Finally, we connect all vertices in $\{x,y\}\cup \{z'_{i1}, z_{i2}', z_{i3}'\,:\, 1\le i\le t\}$
into a cycle $C''$ such that $\{z'_{i1}z_{i2}', z_{i2}'z_{i3}'\,:\, 1\le i\le t\}\subseteq E(C'')$.
If $T'$ is a HIST of $G'$, then $T'':=T'\cup \{z_i'z'_{i1}, z_i'z_{i2}', z_i'z_{i3}'\,:\, 1\le i\le t\}$
is a HIST of $G''$ and $T''\cup C''$ is an SGHG of $G''$.  Conversely,
suppose $H=T\cup C$ is an SGHG of $G''$. We claim that $C=C''$. This in turn gives
that $T=T''$ and therefore $T''-\{z'_{i1}, z_{i2}', z_{i3}'\,:\, 1\le i\le t\}$ is a
HIST of $G'$.  To show that $C=C''$, we first show that
$z_{i2}'\in L(T)$ for each $i$. Suppose on the contrary and assume, without loss of
generality, that $z_{12}'\in S(T)$. Then as $N_{G''}(z_{12}')=\{z_1', z_{11}', z_{13}'\}$,
we get $\{z_{12}'z_1', z_{12}'z_{11}', z_{12}'z_{13}'\}\subseteq E(T)$.
Since $T$ is acyclic, $z_{11}'z_1', z_{13}'z_1'\not\in E(T)$. This in turn
shows that $\{z_1', z_{11}', z_{13}'\}\subseteq L(T)$. However,
$\{z_{12}'z_1', z_{12}'z_{11}', z_{12}'z_{13}'\}$ forms a component
of $T$, showing a contradiction. Then we show that $z_{i1}', z_{i3}'\in L(T)$
for each $i$.
Suppose on the contrary and assume, without loss of
generality, that $z_{11}'\in S(T)$.
By the previous argument,
we have $z_{12}'\in L(T)$. Then $z_1', z_{13}'\in L(T)$
as $z_{12}'$ is on $C$ and $z_1'$ and $z_{13}'$ are
the only two neighbors of $z_{12}'$ which can be  on the cycle $C$.
As $d_{G''}(z_{11}')=3$
and $\{z_{12}', z_1'\}\subseteq N_{G''}(z_{11}')$,
$z_{11}'z_{12}', z_{11}'z_1'\in E(T)$.
Since $z_{12}'\in L(T)$ and $z_1', z_{13}'\in L(T)$,
we get  $z_{12}'z_{13}', z_{12}'z_1', z_1'z_{13}'\not \in E(T)$.
Since $d_{G''}(z_{12}')=d_{G''}(z_{13}')=3$, we have
$z_{12}'z_{13}', z_{12}'z_1', z_1'z_{13}'\in E(C)$.
However, $z_{12}'z_{13}', z_{12}'z_1', z_1'z_{13}'$ forms a
triangle but $|V(C)|\ge 4$, showing a contradiction.
So we have shown that $\{z_{i1}',z_{i2}',z_{i3}'\,:\, 1\le i\le t\}\subseteq L(T)$.
This indicates that in the tree $T-\{z_{i1}',z_{i2}',z_{i3}'\,:\, 1\le i\le t\}$,
each vertex $z_i'$ has degree 1 and no vertices of degree 2.
Hence  $T-\{z_{i1}',z_{i2}',z_{i3}'\,:\, 1\le i\le t\}$ is a HIST of $G'$.

Combining the arguments in the two paragraphs above, we see that
$G$ has a hamiltonian path between $x$ and $y$ if and only if $G''$
has an SGHG. Hence a polynomial SGHG-tester becomes a polynomial
path-tester.
\qqed

Since a generalized Halin graph is 3-connected, the connectivity
requirement in Theorem~\ref{Main-Result2} is necessary. To show that
the minimum degree requirement is best possible, we show
the following proposition.

\begin{PRO}\label{Best-possible-example}
Let $G(A,B)=K_{a,b}$ be a complete bipartite graph with $|A|=a$ and $|B|=b$.
Then $G(A,B)$ has no HIST $T$ with $|L(T)\cap A|=|L(T)\cap B|$ if $b>\frac{3(a-1)}{2}$.
\end{PRO}

If a bipartite graph $G(A,B)$ contains an SGHG $H=T\cup C$,
then  $|L(T)\cap A|=|L(T)\cap B|$.
Thus, by Proposition~\ref{Best-possible-example},
it is easy to see that
the complete bipartite graphs $K_{a,b}$ with $b=\frac{3a-1}{2}$ when $a$ is odd
and $b=\frac{3a-2}{2}$ when $a$ is even does not have an SGHG.
Let $n=a+b$. By direct computation, we get $\delta(K_{a,b})=\frac{2n+1}{5}$ when
$b=\frac{3a-1}{2}$ and $\delta(K_{a,b})=\frac{2n+2}{5}$ when
$b=\frac{3a-2}{2}$. We now prove Proposition~\ref{Best-possible-example}.

\proof[Proof of Proposition~\ref{Best-possible-example}]
Suppose on the contrary that $G(A,B)$ contains a HIST $T$ such that
$|L(T)\cap A|=|L(T)\cap B|$.  Then
\begin{eqnarray*}
  |S(T)\cap B|-|S(T)\cap A| &=& |B|-|L(T)\cap B|-(|A|-|L(T)\cap A)| \\
   &=& |B|-|A|>\frac{3(a-1)}{2}-a=\frac{a-3}{2}.
\end{eqnarray*}
Since $G(A,B)$ is bipartite and $T$ is a HIST of $G(A,B)$, we have
$|S(T)\cap A|\ge 1$. Thus, from the inequalities above, we obtain
$|S(T)\cap B|>(a-1)/2$.
Since $T$ is a HIST, we have $d_T(y)\ge 3$ for each $y\in S(T)\cap B$.
Let $E_{B}=\{e\in E(T)\,:\, e \, \mbox{is incident to a vertex in $S(T)\cap B$}\}$.
Denote by $T'$ the subgraph of $T$ induced on $E_B$.
Notice that $T'$ is a forest of at least $3|S(T)\cap B|$ edges.
Hence $T'$ has at least $3|S(T)\cap B|+1$ vertices. As $T'$
is a bipartite graph with one partite set as $S(T)\cap B$,
and another as a subset of $A$, we conclude that
$|V(T)\cap A|=|V(T)|-|S(T)\cap B|\ge 2|S(T)\cap B|+1$.
Since $|S(T)\cap B|>(a-1)/2$, we then have $|V(T)\cap A|>a$.
This gives a contradiction to the assumption $|A|=a$.
\qqed

%

\section{Notations and Lemmas}

We consider in this paper simple and finite graphs only.
Given a graph $G$, we denote by $V(G)$ and $E(G)$
the vertex set and edge set of $G$, respectively,
and by $e(G)$ the size of $G$.
Let $S\subseteq V(G)$ and $v\in V(G)$. Denote by $G[S]$ the subgraph of $G$
induced on $S$, and denote by  $\Gamma_G(v,S)$  the set of neighbors of $v$ in $S$, and
$deg_G(v,S)=|\Gamma_G(v,S)|$. When $S=V(G)$, we only write $\Gamma_G(v)$ and
$deg_G(v)$. For two subsets $U_1,U_2 \subseteq V(G)$,
let $\delta_G(U_1,U_2)=\min\{deg_G(u_1,U_2)\,:\, u_1\in U_1\}$ and
$\Delta_G(U_1,U_2)=\max\{deg_G(u_1,U_2)\,:\, u_1\in U_1\}$.
Denote by $E_G(U_1,U_2)$  the set of edges with one end in $U_1$ and the other in $U_2$,
the cardinality of $E_G(U_1,U_2)$ is denoted by $e_G(U_1,U_2)$.
Let $u,v\in V(G)$
be two vertices. We write $u\sim v$ if $u$ and $v$ are adjacent.
A path connecting  $u$ and $v$ is called
a  $(u,v)$-path.
If $G$ is a bipartite graph with partite sets
$A$ and $B$, we denote $G$ by $G(A,B)$ for specifying the
two partite sets.
A \emph{matching} in $G$ is a set of independent edges; a \emph{$\wedge$-matching }
is a set of vertex-disjoint copies of $K_{1,2}$; and a \emph{claw-matching}
is a set of vertex-disjoint copies of $K_{1,3}$. The set of degree 2 vertices in
a $\wedge$-matching is called the center of the \emph{$\wedge$-matching }; and
the set of degree 3 vertices in a claw-matching is called the center of the claw-matching.
A cycle $C$ in a graph $G$ is \emph{dominating} if $G-V(C)$ is an edgeless graph.


The Regularity Lemma of Szemer\'edi\,\cite{Szemeredi-regular-partitions} and Blow-up lemma of
Koml\'os et al.\,\cite{Blow-up} are main tools in our proof of Theorem~\ref{Main-Result2}.
For any two disjoint non-empty vertex-sets $A$ and $B$ of a graph $G$, the \emph{density}
of $A$ and $B$ is the ratio $d(A,B):=\frac{e(A,B)}{|A||B|}$. Let $\varepsilon$ and
$\delta$ be two positive real numbers. The pair $(A,B)$ is called $\ve$-regular if for
every $X\subseteq A$ and $Y\subseteq B$ with $|X|>\ve|A|$ and $|Y|>\ve|B|$,
$|d(X,Y)-d(A,B)|<\ve$ holds. In addition, if $deg(a,B)>\delta |B|$ for each $a\in A$ and
$deg(b,A)>\delta|A|$ for each $b\in B$, we say $(A,B)$ an $(\ve,\delta)$-super regular pair.

\begin{LEM}[\textbf{Regularity lemma-Degree form~\cite{Szemeredi-regular-partitions}}]\label{regularity-lemma}
For every $\ve>0$ there is an $M=M(\ve)$ such that if $G$ is any graph with $n$ vertices and $d\in[0,1]$ is
any real number, then there is a partition of the vertex set $V(G)$ into $l+1$ clusters $V_0,V_1,\cdots,V_l$,
and there is a spanning subgraph $G'\subseteq G$ with the following properties.
\vspace{-4mm}
\begin{itemize}
  \item $l\le M$;
  \item $|V_0|\le \ve n$, all clusters $|V_i|=|V_j|\le \lceil\ve n\rceil$ for all $1\le i\ne j\le l$;
  \item $deg_{G'}(v)>deg_G(v)-(d+\ve)n$ for all $v\in V(G)$;
  \item $e(G'[V_i])=0$ for all $i\ge 1$;
  \item all pairs $(V_i,V_j)$ ($1\le i< j\le l$) are $\ve$-regular, each with a density either $0$ or greater than $d$.
 \end{itemize}
\end{LEM}

\begin{LEM}[\textbf{Blow-up lemma-weak version~\cite{Blow-up}}]\label{blow-up}
Given a graph $R$ of order $r$ and positive parameters $\delta, \Delta$, there
exists a positive $\ve=\ve(\delta,\Delta, r)$ such that the following holds.
Let $n_1,n_2,\cdots, n_r$ be arbitrary positive integers and let us replace the
vertices $v_1,v_2,\cdots,v_r$ with pairwise disjoint sets $V_1, V_2,\cdots, V_r$
of sizes $n_1,n_2,\cdots, n_r$\,(blowing up). We construct two graphs on the same vertex set
$V=\bigcup V_i$. The first graph $K$ is obtained by replacing each edge $v_iv_j$ of
$R$ with the complete bipartite graph between the corresponding vertex sets $V_i$ and $V_j$.
A sparser graph $G$ is constructed by replacing each edge $v_iv_j$ arbitrarily with an
$(\ve,\delta)$-super regular pair between $V_i$ and $V_i$. If a graph $H$ with $\Delta(H)\le \Delta$
is embeddable into $K$ then it is already embeddable into $G$.
\end{LEM}

\begin{LEM}[\textbf{Blow-up lemma-strengthened version~\cite{Blow-up}}]\label{blow-up-strengthened}
Given $c>0$, there are positive numbers $\ve=\ve(\delta, \Delta, r,c)$ and $\gamma=\gamma(\delta, \Delta, r,c)$
such that the Blow-up lemma in the equal size case\,(all $|V_i|$ are the same) remains true if for every $i$
there are certain vertices $x$ to be embedded into $V_i$ whose images are a priori restricted to certain sets
$C_x\subseteq V_i$ provided that
\vspace{-5mm}
\begin{enumerate}
  \item [(i)] each $C_x$ within a $V_i$ is of size at least $c|V_i|$;\vspace{-3mm}
  \item [(ii)] the number of such restrictions within a $V_i$ is not more than $\gamma|V_i|$.
 \end{enumerate}
 \end{LEM}
\vspace{-4mm}
We will use both the weak and strengthened versions of Blow-up lemma in our proof.

Besides the above two lemmas, we also need the two lemmas below regarding regular pairs.

\begin{LEM}\label{regular-pair-large-degree}
If $(A,B)$ is an $\ve$-regular pair with density $d$, then for any $A'\subseteq A$ with
$|A'|> \ve |A|$, there are at most $\ve |B|$ vertices $b\in B$ such that
$deg(b, A')\le (d-\ve)|A'|$.
\end{LEM}

\begin{LEM}[\textbf{Slicing lemma}]\label{slicing lemma}
Let $(A,B)$ be an $\ve$-regular pair with density $d$, and for some $\nu >\ve$, let $A'\subseteq A$
and $B'\subseteq B$ with $|A'|\ge \nu|A|$, $|B'|\ge \nu|B|$. Then $(A',B')$ is an $\ve'$-regular
pair of density $d'$, where $\ve'=\max\{\ve/\nu, 2\ve\}$ and $d'>d-\ve$.
\end{LEM}


The following two results on hamiltonicity are used in finding cycles in the proofs.

\begin{LEM}[\cite{Ore-Hamiltonian}]\label{Hamiltonian-connected}
If G is a graph of order $n$ satisfying $d(x) + d(y)\ge  n + 1$ for every pair of
nonadjacent vertices $x, y \in  V(G)$, then G is hamiltonian-connected.
\end{LEM}

\begin{LEM}[\cite{M-M-bipartite-hamiltonian}]\label{bi-Hamiltonian}
Let $G$ be a balanced bipartite graph with $2n$ vertices. If $d(x)+d(v)\ge n+1$
for any two non-adjacent vertices $x,y\in V(G)$, then $G$ is hamiltonian.
\end{LEM}

\section{Proof of Theorem~\ref{Main-Result2}}

Given $0\le \beta\ll\alpha\ll 1$, we define the two extremal cases with parameters $\alpha$
and $\beta$
as follows.

\noindent
\textbf{Extremal Case 1.} There exists a partition of $V(G)$ into $V_1$ and $V_2$
such that $|V_i|\ge(2/5-4\beta)n$ and
$d(V_1,V_2)< \alpha$. Furthermore, $deg(v_1,V_2)\le  2\beta n$ for each $v_1\in V_1$.

\noindent
\textbf{Extremal Case 2.} There exists a partition of $V(G)$ into $V_1$ and $V_2$
such that $|V_1|>(3/5-\alpha)n$ and $d(V_1,V_2)\ge 1-3\alpha$. Furthermore,
 $deg(v_1,V_2)\ge (2n+3)/5-2\beta n$ for each $v_1\in V_1$.

Then Theorem~\ref{Main-Result2} is shown through the following three theorems.

\begin{THM}[\textbf{Non-extremal Case}]\label{non-extremal}
For every $\alpha>0$, there exists $\beta>0$ and a positive integer $n_0$ such that
if $G$ is a 3-connected graph with $n\ge n_0$ vertices and $\delta(G)\ge (2n+3)/5-\beta n$, then
$G$ contains an SGHG or $G$ is in one of the two extremal cases.
\end{THM}

\begin{THM}[\textbf{Extremal Case 1}]\label{extremal_1}
Suppose that $0<\beta \ll \alpha\ll 1$ and $n$ is a sufficiently large integer.
Let $G$ be a 3-connected graph on $n$ vertices with $\delta(G)\ge (2n+3)/5$.
If $G$ is in Extremal Case 1, then $G$ contains an SGHG.
\end{THM}

\begin{THM}[\textbf{Extremal Case 2}]\label{extremal_2}
Suppose that $0<\beta \ll \alpha\ll 1$ and $n$ is a sufficiently large integer.
Let $G$ be a 3-connected graph on $n$ vertices with $\delta(G)\ge (2n+3)/5$.
If $G$ is in Extremal Case 2, then $G$ contains an SGHG.
\end{THM}

We show Theorems~\ref{non-extremal}-\ref{extremal_2} separately in the
following three subsections.

\subsection{Proof of Theorem~\ref{non-extremal}}
We fix the following sequence of parameters,
\begin{equation}\label{parameters}
   0<\ve\ll d\ll \beta\ll \alpha < 1,
\end{equation}
and specify their dependence as the proof proceeds. We let $\beta\ll\alpha$ be the same $\alpha$
and $\beta$ as
defined in the two extremal cases. Then we choose $d\ll \beta$. Finally we choose
\begin{equation}\label{epsilon}
   \ve=\min\left\{\frac{1}{4}\ve\left(\frac{d}{2}, \left\lceil\frac{2}{d^3}\right\rceil, 2,\frac{d}{2}\right), \frac{1}{9}\ve\left(\frac{d}{2}, \left\lceil\frac{3}{d^3}\right\rceil, 3\right), \frac{1}{4}\ve\left(\frac{d}{2}, 2, 2,\frac{d}{2}\right)  \right \},
\end{equation}
where $\ve\left(\frac{d}{2}, \left\lceil\frac{3}{d^3}\right\rceil, 3\right)$ follows
from the definition of the $\ve$ in the weak version of  the Blow-up lemma
and $\ve\left(\frac{d}{2}, \left\lceil\frac{2}{d^3}\right\rceil, 2,\frac{d}{2}\right)$ and
$\ve\left(\frac{d}{2}, 2, 2,\frac{d}{2}\right)$
follow
from the definition of the $\ve$ in the  strengthened version of the Blow-up lemma.
Choose $n$ to be sufficiently large. In the proof, we omit non-necessary ceiling and floor functions.

Let $G$ be a graph of order $n$ such that $\delta(G)\ge (2n+3)/5-\beta n$ and suppose that
$G$ is not in any of the
two extremal cases. Applying the regularity lemma to $G$ with parameters $\ve$ and $d$, we obtain a partition of $V(G)$
into $l+1$ clusters $V_0,V_1,\cdots, V_l$ for some $l\le M=M(\ve)$, and a spanning subgraph $G'$ of $G$ with all described
properties in Lemma~\ref{regularity-lemma}\,(the Regularity lemma). In particular, for all $v\in V$,
\begin{eqnarray}\label{G'degree}
  deg_{G'}(v) &>& deg_G(v)-(d+\ve)n\ge(2/5-\beta-d-\ve)n \nonumber \\
   &\ge& (2/5-2\beta)n\,\,(\mbox{provided that }\ve+d\le \beta),
\end{eqnarray}
and
$$
e(G')\ge e(G)-\frac{(d+\ve)}{2}n^2\ge e(G)-dn^2,
$$
by  using $\ve<d$.

We further assume that $l=2k$ is even; otherwise, we eliminate the last cluster $V_l$ by removing all
the vertices in this cluster to $V_0$. As a result, $|V_0|\le 2\ve n$ and
\begin{equation}\label{lN and n}
   (1-2\ve)n\le lN=2kN\le n,
\end{equation}
here we assume that $|V_i|=N$ for $i\ge 1$.

For each pair $i$ and $j$ with $1\le i< j\le l$, we write $V_i\sim V_j$
 if $d(V_i,V_j)\ge d$. We now consider the reduced graph $G_r$, whose vertex set is $\{1,2,\cdots, l\}$, and
 two vertices $i$ and $j$ are adjacent if and only if $V_i\sim V_j$. We claim that $\delta(G_r)\ge (2/5-2\beta)l$.
 Suppose not, and let $i_0\in V(G_r)$ such that  $deg(i_0,V(G_r))<(2/5-2\beta)l$. Then, for the corresponding cluster
 $V_{i_0}$ we have $e_{G'}(V_{i_0}, V(G')-V_{i_0}) < |V_{i_0}|(2/5-2\beta)lN$.
 On the other hand, by using (\ref{G'degree}), we have  $e_{G'}(V_{i_0}, V(G')-V_{i_0})\ge |V_{i_0}|(2/5-2\beta)n$.
 As $lN\le n$ from~(\ref{lN and n}), we obtain a contradiction.
The rest of the proof consists of the following steps.

\noindent \textbf{Step 1.} Show that $G_r$ contains a dominating cycle $C$ and there is a
$\wedge$-matching in $G_r$ with all vertices in $V(G_r)-V(C)$ as its center.
We distinguish two cases in Step 1, and each of the  other steps  will be separated into two cases correspondingly.

\noindent \textbf{Case A.} $C=X_1Y_1X_2Y_2\cdots X_tY_t$ is an even cycle for some $t\le k$.

\noindent \textbf{Case B.} $C=X_0X_1Y_1X_2Y_2\cdots X_tY_t$ is an odd cycle for some $t< k$.

Notice that in Case B there is at least one vertex in $V(G_r)-V(C)$ by the assumption that $|V(G_r)|=l$
is even. In what follows, if we denote a vertex of $G_r$ by a capital letter, it means either
a vertex of $G_r$ or the corresponding cluster in $G$, but the exact meaning
will be clear from the context.
For $1\le i\le t$, we call $X_i$ and $Y_i$ the partners of each other, and write as $P(X_i)=Y_i$ and  $P(Y_i)=X_i$.

Since $C$ is not necessarily hamiltonian in $G_r$, we need to take care of the clusters
of $G$ which are not represented on $C$.  For each vertex $F\in V(G_r)-V(C)$, we partition
the corresponding cluster $F$ into two small clusters $F_1$ and $F_2$
such that $-1 \le |F_1|-|F_2|\le 1$. We call each  $F_1$ and $F_2$
a {\it half-cluster}.  Then we group all the original clusters and
the partitioned clusters into pairs $(A,B)$ and triples $(C,D,F)$ with $F$
as a half-cluster
such that each pair  $(A,B)$ and $(C,D)$ is still $\ve$-regular with density $d$
and the pair $(D,F)$ is $2.1\ve$-regular with density $d-\ve$. Having the
cluster groups like this, in the end, we will find ``small'' HITs
within each pair $(A,B)$ or among each triple $(C,D,F)$.

\noindent \textbf{Step 2.}
For each $1\le i\le t-1$, initiate two independent edges
connecting $Y_{i}$ and $X_{i+1}$. In Case A, also
initiate two independent edges
connecting $X_1$ and $Y_t$; and in Case B,
initiate two independent edges
connecting  the clusters in
each pair of $X_0$ and $X_1$, and $X_0$ and $Y_t$.


\noindent \textbf{Step 3.} Make each regular pair in the new grouped pairs and triples given in Step 1 super-regular.

\noindent \textbf{Step 4.} Construct HITs covering all vertices in $V_0$ using vertices from the super-regular pairs obtained
from Step 3, and obtain new  super-regular pairs.

\noindent \textbf{Step 5.} Apply the Blow-up lemma to find a HIT between a super-regular pair resulted from Step 4 or among a triple $(A,B,F)$,
where both $(A,F)$ and $(A,B)$ are super-regular pairs resulted from Step 4, and $F$
is a half cluster. In addition, in the construction, for each triple $(A,B,F)$,
we require the HIT to use as many vertices as possible from $F$ as non-leaves.

\noindent \textbf{Step 6.}  Apply the Blow-up Lemma again on the regular-pairs induced on the leaves of each HIT obtained in Step 5
to find two disjoint paths covering all the leaves. Then connect all the HITs into a HIST of $G$ using edges
guaranteed by the regularity and connect the
disjoint paths into a cycle  using the edges initiated in Step 2. The union of the HIST and the cycle gives an SGHG of $G$.
%

We now give details of each step.
The assumption that $G$ is not in any of the two extremal cases leads to the following claim, which will be used in Step 1.

\begin{CLA}\label{extremal-claim}
Each of the following holds for $G_r$.
\vspace{-3mm}
\begin{itemize}
  \item [$($a$)$] $G_r$ contains no cut-vertex set of size at most $\beta l$; \vspace{-3mm}
  \item [$($b$)$]  $G_r$ contains no independent set  of size more than $(3/5-\alpha/2)l$ .
\end{itemize}
\end{CLA}
\pf (a) Suppose instead that $G_r$ contains a vertex-cut $W$ of size at most $\beta l$. As $\delta(G_r)\ge (2/5-2\beta)l$,
then each component of $G_r-W$ has at least $(2/5-3\beta)l$ vertices. Let $U$ be the vertex set of one of the components
of $G_r-W$, $A=\bigcup_{i\in U}V_i$, and $B= V(G)-A$.
We see that $|A|, |B|\ge (2/5-3\beta)lN\ge (2/5-4\beta)n$, and since $e(G)\le e(G')+dn^2$, we have
\begin{eqnarray*}
  e_G(A,B) &\le & e_{G'}(A,B)+dn^2\le |W||A|+dn^2 \\
   &\le & \beta l N(3/5+3\beta)lN+dn^2\le (3\beta/5+3\beta^2+d)n^2 \quad (\mbox{as}\,|A|\le (3/5+3\beta)lN\,\mbox{and}\, ln\le n)\\
   & \le & \frac{25}{3}(3\beta/5+3\beta^2+d)|A||B|\quad (\mbox{since}\,|A||B|\ge 3n^2/25)\\
   &<& \alpha |A||B|\quad (\mbox{provided that} \, \frac{25}{3}(3\beta/5+3\beta^2+d)<\alpha).
\end{eqnarray*}
This shows that $d(A,B)<\alpha$. Since $deg_{G_r}(u, V(G_r)-U)=deg_{G_r}(u,W)\le \beta l $
for each $u\in U$,
we see that $deg_G(a, B)\le \beta lN+(d+\ve)n\le 2\beta n$  for each $a\in A$ provided that $d+\ve \le \beta$.
However, the above argument shows that $G$ is in Extremal Case 1, showing a contradiction.

(b) Suppose instead that $G_r$ contains an independent set $U$ of size larger than $(3/5-\alpha/2)l$.
Let $U'=V(G_r)-U$, $A=\bigcup_{i\in U}V_i$, and $ B= V(G)-A$.
Then $|A|\ge (3/5-\alpha/2)lN\ge (3/5-\alpha)n$.
For each vertex $v\in A$, since $deg_G(v,A)\le deg_{G'}(v,A)+(d+\ve)n\le
\beta n$, we have
$deg_G(v,B)\ge (2n+3)/5-\beta n-\beta n \ge (2n+3)/5-2\beta n$. This gives that
$$
d(A,B)\ge \frac{(2/5-2\beta)n}{|B|}\ge \frac{(2/5-2\beta)n}{(2/5+\alpha)n}\ge 1-3\alpha,
$$
provided that $\beta \le \alpha/10+3\alpha^2/2$.
We see that $G$ is in Extremal Case 2.
\qed

\noindent \textbf{Step 1.} Show that $G_r$ contains a dominating cycle $C$, and there is a
$\wedge$-matching in $G_r$ with all vertices in $V(G_r)-V(C)$ as its center.

We need some results on longest cycles and paths as follows.

\begin{LEM}[\cite{Nash-Williams-dominating}]\label{2/3n-dominating-cycle}
Let $G$ be a 2-connected graph on $n$ vertices with $\delta(G)\ge (n+2)/3$. Then every longest cycle in
$G$ is a dominating cycle.
\end{LEM}

\begin{LEM}[\cite{Bauer-long-cycles}]\label{Bauer-long-cycles}
Let $G$ be a 2-connected graph on $n$ vertices with $\delta(G)\ge (n+2)/3$. Then $G$ contains a cycle of length
at least $\min\{n,n+\delta(G)-\alpha(G)\}$, where $\alpha(G)$ is the size of a largest independent set in $G$.
\end{LEM}

\begin{LEM}[\cite{Relative-length-diff}]\label{diff}
If $G$ is a 3-connected graph of order $n$ such that the degree sum of any four independent vertices is
at least 3n/2+1, then the number of vertices on a
longest path and that on a longest cycle differs at most by 1.
\end{LEM}

By (a) of Claim~\ref{extremal-claim},
$G_r$ is $\beta l$-connected. Since $n=Nl+|V_0|\le (l+2)\ve n$,  we get $l\ge 1/\ve -2$.
Since $1/\ve-2\ge 3/\beta$\,(provided that $\beta \ge 3\ve/(1-2\ve)$), we then have
$\beta l\ge 3$.  So  $G_r$ is 3-connected. By Claim~\ref{extremal-claim} (b), $G_r$ has no independent set of size
more  than $(3/5-\alpha/2)l$. Notice that $\delta(G_r)\ge (2/5-2\beta)l>(l+2)/3$.
Applying Lemma~\ref{2/3n-dominating-cycle} and Lemma~\ref{Bauer-long-cycles}
on $G_r$, we see that  there is a cycle $C$ in $G_r$ which is longest, dominating, and has length
at least $(4/5+\alpha/2-2\beta)l $. Let $\mathcal{W}=V(G_r)-V(C)$.
In Case B, we order and label the vertices of $C$ such that $X_0$ is adjacent
to a vertex, say $Y_0\in \mathcal{W}$\,(recall that $\mathcal{W}\ne \emptyset$ in this case).
We fix $(X_0,Y_0)$
as a pair at the first place\,($X_0Y_0\in E(G_r)$, as cluster in $G$, $(X_0,Y_0)$ is an $\ve$-regular pair with density $d$). Let
$$
\mathcal{W}'=
\left\{
  \begin{array}{ll}
    \mathcal{W}, & \hbox{if in Case A;} \\
    \mathcal{W}-\{Y_0\}, & \hbox{if in Case B.}
  \end{array}
\right.
$$
We have
$|\mathcal{W}'|\le (1/5-\alpha/2+2\beta)l$ if in Case A and $|\mathcal{W'}|\le (1/5-\alpha/2+2\beta)l-1$
if in Case B. So
 $2|\mathcal{W}'|\le (2/5-\alpha+4\beta)l<(2/5-2\beta)l $\,(provided that $\beta<\alpha/6$) if in Case A
 and $2|\mathcal{W}'|\le (2/5-\alpha+4\beta)l-2<(2/5-2\beta)l-1 $\,(provided that $\beta<\alpha/6$) if in Case B.
Thus
there is a $\wedge$-matching centered in all vertices in $\mathcal{W'}$; furthermore,
if in Case B, we can choose the matching such that
 $X_0$ is not covered by it. Let $M_{\wedge}$ be  such a matching.
For a vertex $X\in \mathcal{W}'$, denote by $M_{\wedge}(X)$  the two vertices from $V(C)$ to which $X$ is adjacent in
$M_{\wedge}$. Then we have two facts about vertices in $M_{\wedge}(X)$.

\begin{FAC}\label{no-2-consecutive}
Let $X\in \mathcal{W}'$. Then the two vertices in $M_{\wedge}(X)$ are non-consecutive on $C$.\,
$($By the assumption that $C$ is longest.$)$
\end{FAC}

\begin{FAC}\label{no-adjacent-dominator}
Let $X\ne Y\in \mathcal{W}'$. Then no two vertices from $M_{\wedge}(X)\cup M_{\wedge}(Y)$
are adjacent  on $C$. $($By applying Lemma~\ref{diff}.$)$
\end{FAC}
%
%
%
%

For a complete bipartite graph, if it contains an SGHG,
then the ratio of the cardinalities of the two partite sets
should be greater than $2/3$ as shown in Proposition~\ref{Best-possible-example}.
Since a longest dominating cycle in $G_r$ is not necessarily hamiltonian,
we need to take care of the clusters of $G$ which are not represented
by the vertices on $C$. One possible consideration is that
for each $F\in V(G_r)-V(C)$, suppose $F$ is adjacent to
$A\in V(C)$, recall $P(A)$ is the partner of $A$. Then as clusters, we consider the bipartite
graph of $G$ with partite sets $A$ and $P(A)\cup F$.
However,  $|A|/|P(A)\cup F|$  is about 1/2, which is less
than $2/3$.
For this reason,  we partition $F\in V(G_r)-V(C)$ into two parts to attain the right ratio
in the corresponding bipartite graphs.
Suppose $M_{\wedge}(F)=\{D_1,D_2\}\subseteq V(C)$. As a cluster of
$G$, we partition $F$ into $F_1$ and $F_2$ arbitrarily such that
$$
|F_1|=\left\lfloor\frac{|F|}{2}\right\rfloor=\left\lfloor\frac{N}{2}\right\rfloor \quad \mbox{and} \quad
|F_2|=\left\lceil\frac{|F|}{2}\right\rceil=\left\lceil\frac{N}{2}\right\rceil.
$$
We call each $F_i$ a half-cluster of $G$.  Then we create two pairs $(D_i,F_i)$, and call $D_i$ the dominator
of $F_i$, and $F_i$ the follower of $D_i$, and $(D_i,F_i)$ a DF-pair, for $i=1,2$.
We have the following fact
about a DF-pair.
\begin{FAC}\label{DF-pair-1level}
Each DF-pair $(D,F)$ is $2.1\ve$-regular with density at least $d-\ve$. $($By Slicing lemma.$)$
\end{FAC}

Also, by Fact~\ref{no-2-consecutive} and Fact~\ref{no-adjacent-dominator},
if $D\in V(C)$  is a dominator, then $P(D)$, the partner of $D$, is not a dominator for any followers.
As $X_0\not\in V(\mathcal{W}')$, we know that $X_0$ is not
a dominator for any half-clusters.
We group the clusters and half-clusters of $G$ into {\it H-pairs} and {\it H-triples} in a way below.
For each pair  $(X_i,Y_i)$ on $C$, if $\{X_i,Y_i\}\cap V(M_{\wedge})=\emptyset$, we take $(X_i,Y_i)$ as an H-pair.
Otherwise, $|\{X_i,Y_i\}\cap V(M_{\wedge})|=1$ by Fact~\ref{no-2-consecutive} and Fact~\ref{no-adjacent-dominator}.
Since there is no difference for the proof for the case that $X_i\in V(M_{\wedge}) $  or the case that
$Y_i\in V(M_{\wedge})$, throughout the remaining proof, we always  assume that  $Y_i\in V(M_{\wedge})$
if $\{X_i,Y_i\}\cap V(M_{\wedge})\ne\emptyset$. In this case, there is a unique half-cluster $F$
with $Y_i$ as its dominator. Then we take $(X_i,Y_i, F)$ as an H-triple. We assign $(X_0,Y_0)$ as an $H$-pair.

\noindent \textbf{Step 2.} Initiating connecting edges.


Given an $\ve$-regular pair $(A,B)$ of density $d$ and a subset $B'\subseteq B$, we say a vertex $a\in A$
\emph{typical} to $B'$ if $deg(a, B')\ge (d-\ve)|B'|$. Then by the regularity of $(A,B)$, the fact below
holds.

\begin{FAC}\label{Typical}
If $(A,B)$ is an $\ve$-regular pair, then at most $\ve|A|$ vertices of $A$ are not typical to $B'\subseteq B$
whenever $|B'|>\ve |B|$.
\end{FAC}


For each $1\le i\le t-1$,
choose $y^*_{i}\in Y_{i}$ typical to both $X_i$
and $X_{i+1}$, and
$y^{**}_{i}\in Y_{i}$ typical to each of $X_i$, $X_{i+1}$,  and $\Gamma(y_i^*, X_i)$.
Correspondingly, choose $x_{i+1}^*\in \Gamma(y^*_{i}, X_{i+1})$
typical to $Y_{i+1}$, and $x_{i+1}^{**}\in \Gamma(y^{**}_{i}, X_{i+1})$
typical to both $Y_{i+1}$ and $\Gamma(x_{i+1}^*, Y_{i+1})$.
For $i=t$, we choose $y_t^*$ and $y_t^{**}$
the same way as for $i<t$,
but if in Case A,
choose $x_{1}^*\in \Gamma(y^{**}_{t}, X_{1})$
typical to $Y_{1}$, and $x_{1}^{**}\in \Gamma(y^{*}_{t}, X_{1})$
typical to both $Y_{1}$ and $\Gamma(x_{1}^*, Y_{1})$;
and if in Case B,
choose $x_{0}^*\in \Gamma(y^{**}_{t}, X_{0})$
typical to $X_{1}$, and $x_{0}^{**}\in \Gamma(y^{*}_{t}, X_{0})$
typical to both $X_{1}$ and $\Gamma(x_{0}^*, X_{1})$.
Furthermore, in Case B,  we  choose $y^*_{t+1}\in X_{0}$ typical to both $Y_0$
and $X_{1}$, and
$y^{**}_{t+1}\in X_{0}$ typical to each of $Y_0$, $X_{1}$,  and $\Gamma(y_{t+1}^*, Y_0)$.
Correspondingly, choose $x_{1}^*\in \Gamma(y^*_{t+1}, X_{1})$
typical to $Y_1$
and  $x_{1}^{**}\in \Gamma(y^{**}_{t+1}, X_{1})$ typical to both $Y_1$ and $\Gamma(x_1^*, Y_1)$.
Additionally, we choose $y_0^*\in \Gamma(y_{t+1}^*, Y_0)$ such that $y_0^*$ is typical to
$X_0$, and choose $y_0^{**}\in  \Gamma(y_{t+1}^{**}, Y_0)$ such that $y_0^{**}$ is typical to
$X_0$. Notice that by the choice of these vertices above, we have the following.
$$
\left\{
  \begin{array}{ll}
    y_i^*x_{i+1}^*, y_i^{**}x_{i+1}^{**}\in E(G), & \hbox{for $1\le i\le t-1$;} \\
    x_1^*y_t^{**}, x_1^{**}y_t^*\in E(G), & \hbox{in Case A;} \\
    x_0^*y_t^{**}, x_0^{**}y_t^*, x_1^*y_{t+1}^*, x_1^{**}y_{t+1}^{**}, y_0^*y_{t+1}^*, y_0^{**}y_{t+1}^{**}\in E(G), & \hbox{in Case B.}
  \end{array}
\right.
$$

By Fact~\ref{Typical}, for each $0\le i\le t$,
we have
$|\Gamma(x_i^*, Y_i)\cap \Gamma(x_i^{**}, Y_i)|, |\Gamma(y_i^*, X_i)\cap \Gamma(y_i^{**}, X_i)|\ge (d-\ve)^2N$,
and $|\Gamma(y_{t+1}^*, Y_0)\cap \Gamma(y_{t+1}^{**}, Y_0)|\ge (d-\ve)^2N$.



\noindent \textbf{Step 3.} Super-regularizing the regular pairs in each H-pair and H-triple  given in Step 1.

For each $0\le i\le t$, if $(X_i,Y_i)$ is an $H$-pair, let
$$
  X_i'= \{x\in X_i: deg(x,Y_i)\ge (d-\ve)N\}\quad  \mbox{and}\quad Y_i' = \{y\in Y_i: deg(y,X_i)\ge (d-\ve)N\}. \\
$$
By Fact~\ref{Typical}, we have
 $|X_i'|, |Y_i'|\ge (1-\ve)N$.
 Recall that $x_i^*, x_i^{**}\in X_i$ and $y_i^*, y_i^{**}\in Y_i$
are the initiated vertices in Step 2.
 For $1\le i\le t $,  if
 $|X'_i-\{x_i^*,x_i^{**}\}|\ne |Y'_i-\{y_i^*, y_i^{**}\}|$, say
$|X'_i-\{x_i^*,x_i^{**}\}|>|Y'_i-\{y_i^*, y_i^{**}\}|$,
we then remove $|X'_i-\{x_i^*,x_i^{**}\}|-|Y'_i-\{y_i^*, y_i^{**}\}|$
 vertices out from $X_i'-\{x_i^*,x_i^{**}\}$, and denote the remaining set
 still as $X_i'$. Denote $Y'_i-\{y_i^*, y_i^{**}\}$
 still as $Y_i'$. We see that $|X'_i|=|Y'_i|$.
 As $|Y_i'|\ge (1-\ve)N$\,(to be precise, the lower bound should be
 $(1-\ve)N-2$, however, the constant 2 can be made vanished by adjusting the $\ve$ factor,
 we ignore the slight different of the $\ve$-factor here),
  we have that $|X_i\cup Y_i-(X_i'\cup Y_i')|\le 2\ve N$.
 For $i=0 $, if $|X'_i-\{x_i^*,x_i^{**}, y_{t+1}^*, y_{t+1}^{**}\}|\ne |Y'_i-\{y_i^*, y_i^{**}\}|$,
say
$|X'_i-\{x_i^*,x_i^{**}, y_{t+1}^*, y_{t+1}^{**}\}|>|Y'_i-\{y_i^*, y_i^{**}\}|$,
then we remove $|X'_i-\{x_i^*,x_i^{**}, y_{t+1}^*, y_{t+1}^{**}\}|-|Y'_i-\{y_i^*, y_i^{**}\}|$
 vertices out from $X_i'-\{x_i^*,x_i^{**}, y_{t+1}^*, y_{t+1}^{**}\}$ and denote the remaining set
 still as $X_i'$. Denote $Y'_i-\{y_i^*, y_i^{**}\}$
 still as $Y_i'$. We see that $|X'_i|=|Y'_i|$.
We call the resulting  H-pairs {\it supper-regularized H-pairs}.
By Slicing lemma\,(Lemma~\ref{slicing lemma}) and the definitions of $X_i',Y_i'$,
we see that
\begin{FAC}\label{super-regularized XiYi}
Each supper-regularized H-pair $(X_i',Y_i')$ is a $(2\ve, d-2\ve)$-super-regular pair.
\end{FAC}
For each H-triple $(X_i,Y_i, F)$, by Fact \ref{DF-pair-1level},  $(Y_i,F)$ is $2.1\ve$-regular with density at least $d-\ve$.
Let
\begin{eqnarray*}
  X_i' &=&\{x\in X_i:  deg(x,Y_i)\ge (d-\ve)N\}, \\
  Y_i' &=& \{y\in Y_i: deg(y,X_i)\ge (d-\ve)N, deg(y,F)\ge (d-3.1\ve)|F|\},\mbox{and} \\
  F' &=& \{f\in F: deg(f,Y_i)\ge (d-3.1\ve)N\}.
\end{eqnarray*}
Recall that $x_i^*, x_i^{**}\in X_i$ and $y_i^*, y_i^{**}\in Y_i$
are the initiated vertices in Step 2. We remove $x_i^*, x_i^{**}$ out from $X_i'$,
and remove $y_i^*, y_i^{**}$ out from $Y_i'$. Still denote the resulted
clusters as $X_i'$ and $Y_i'$, respectively.
Remove $\lceil d^3N\rceil$ vertices out from $F$, which consists of all vertices in $F-F'$ and
any $\lceil d^3N\rceil-|F-F'|$ vertices from $F'$\,(we need to increase the
ratio $|Y_i'|/|X_i'\cup F'|$ a little as later on we may use vertices in $Y_i'$ in constructing HITs covering vertices in
 $V_0$).  Denote the resulting set still by $F'$.
Then we see that $|X_i'|\ge (1-\ve)N$, $|Y_i'|\ge (1-3.1\ve)N$,
and $|F'|\ge (1-2.1\ve)|F|-d^3N\ge (1-2.1\ve-2d^3)|F|$.
We call the resulted H-triples {\it supper-regularized H-triples}.
By the Slicing Lemma and the definitions above, the following is true.
\begin{FAC}\label{super-regularized H-triple}
For each super-regularized H-triple
$(X_i',Y_i', F')$,
$(X_i', Y_i')$ is
$(2\ve, d-3.1\ve)$-
super-regular,
and $(Y_i',F')$ is $(4.2\ve, d-3.1\ve-2d^3)$-super-regular.
\end{FAC}


Let $V^1_0$ be the union of the set of vertices from each
$(X_i\cup Y_i-(X_i'\cup Y_i'))-\{x_i^*,x_i^{**}, y_i^*, y_i^{**}\}-\{y_{t+1}^*, y_{t+1}^{**}\}$
\,($\{y_{t+1}^*, y_{t+1}^{**}\}$ exists only if in Case B),
where $(X_i, Y_i)$ is an H-pair,
and let $V^2_0$ be the union of the  set of vertices from each
$(X_i\cup Y_i\cup F-(X_i'\cup Y_i'\cup F'))-\{x_i^*,x_i^{**}, y_i^*, y_i^{**}\}$,
where $(X_i,Y_i,F)$ is an H-triple.
Notice that for each H-pair $(X_i,Y_i)$, we have $|X_i\cup Y_i-(X_i'\cup Y_i')|\le 2\ve N$;
and for each H-triple $(X_i,Y_i, F)$, we have $|X_i-X_i'|\le \ve N$, $|Y_i-Y_i'|\le (\ve +2.1\ve) N$, and $|F-F'|\le d^3N$.
Hence by using the facts that
$|\mathcal{W'}|\le (1/5-\alpha/2+2\beta)l$, $t=l/2$, and $Nl\le n$ from
inequality~(\ref{lN and n}), we get
$$
 |V^1_0|+|V^2_0|\le 2\ve N l/2+ 2(1/5-\alpha+2\beta)l(d^3N+  2.1\ve N)\le 2d^3Nl/5+2\ve Nl\le 2d^3n/5+2\ve n.
$$
Let $V_0'=V_0\cup V_0^1\cup V_0^2$. Then
\begin{equation}\label{size V0prime}
 |V'_0|\le 2\ve n+ 2d^3n/5+2\ve n\le  d^3n/2\quad (\mbox{provided that $\ve \le  d^3/40$}).
\end{equation}

\noindent \textbf{Step 4.} Construct small HITs covering all vertices in $V'_0$.

Consider a vertex $x\in V_0'$ and a cluster or a half-cluster $A$, we say that
$x$ is \emph{adjacent to} $A$, denoted by $x\sim A$, if  $deg(x,A)\ge (d-\ve)|A|$.  We call $A$
the \emph{partner} of $x$.

\begin{CLA}\label{absorbing}
For each vertex $x\in V_0'$, there is a cluster or a half-cluster $A$ such that
$x\sim A$, where $A$ is not a dominator, and  we can
assign all vertices in $V_0'$ to their partners which are not dominators
 such that each of the cluster or half-cluster is
used by at most $\frac{d^2N}{20}$ vertices from $V_0'$.
\end{CLA}

\pf Suppose we have found partners for the first $m<d^3n/2$\,(recall that $|V_0'|\le d^3n/2$) vertices
of $V_0'$ such that no cluster or half-cluster is used by at most  $\frac{d^2N}{20}$ vertices.
Let $\Omega$ be the set of all clusters and half-clusters that are used exactly by $\frac{d^2N}{20}$ vertices. Then
\begin{eqnarray*}
  \frac{d^2N}{20}|\Omega| &\le & m <d^3n/2 \le d^3(2kN+2\ve n)/2 \\
   &\le & d^3kN + d^3\frac{2kN}{1-2\ve},
 \end{eqnarray*}
by inequality~(\ref{lN and n}).
Therefore,
\begin{eqnarray*}
  |\Omega| &\le &\frac{20d^3k}{d^2} + \frac{20d^3l}{d^2(1-2\ve)} \\
   &\le & 10dl+ 40 dl \,\,(\mbox{provided that $1-2\ve \ge 1/2$ })\\
   & \le & \beta l \,\,(\mbox{provided that $50d\le \beta$ }).
 \end{eqnarray*}
Consider now a vertex $v\in V_0'$ not having a partner found so far. Let $\mathcal{U}$
be the set of all non-dominator clusters and half-clusters  adjacent to $v$
not contained in $\Omega$. We claim that $|\mathcal{U}|\ge (\alpha-7\beta)l$.
To see this, we first observe that any vertex $v\in V_0'$ is adjacent to
at least $(\alpha-6\beta)l$ non-dominator clusters and half-clusters. For instead,
as $v$ may adjacent to $2|\mathcal{W'}|$ dominators, vertices in $V_0'$,
or clusters $A$ with $deg(v,A)< (d-\ve)|A|$, we have
\begin{eqnarray*}
  (2/5-\beta)n &\le &deg_G(v) < (\alpha-6\beta)lN+ (2/5+4\beta-\alpha)lN+ d^3n/2+(d-\ve)lN \\
   &\le& (2/5-2\beta+d^3/2+d-\ve)n \\
   &<&  (2/5-3\beta/2)n \,\,(\mbox{provided that $d-\ve+ d^3/2< \beta/2$ } ),
\end{eqnarray*}
showing a contradiction. Since $|\Omega|\le \beta l$, we have $|\mathcal{U}|\ge (2\alpha-7\beta)l$.
\qed

Now for each non-dominator cluster $A$\,($A$ is either a cluster $X_i'$, $Y_i'$, or a half cluster $F'$),
let $I(A)$  be the set of vertices from $V_0'$ such that each of them has $A$ as its partner.
By Claim~\ref{absorbing}, we have $|I(A)|\le \frac{d^2N}{20}$.

We need three operations below for constructing small HITs covering vertices in $V_0'$.

\noindent \textbf{Operation I}\quad
Let $(A,B)$ be an $(\ve',\delta)$-super-regular pair, and $I$ a set of vertices disjoint from $A\cup B$. Suppose
that  (i) $ deg(x,B)\ge d'|B|> \ve'|B|$ and $ deg(x,B)\ge d'|B|\ge 3|I|$  for any $x\in I$;
(ii) $(\delta-\ve')d'|B|\ge 3|I|$; (iii) $(\delta-\ve')|A|>|I|$; and
(iv) $\delta|A|>4|I|$.
Then we can do the following operations on $(A,B)$ and $I$.

Let $I=\{x_1,x_2,\cdots, x_{|I|}\}$.
We first assume that $|I|\ge 2$.

\begin{figure}[!htb]
\psfrag{B}{$B$} \psfrag{A}{$A$}
\psfrag{|I|=3}{$|I|=3$} \psfrag{|I|=1}{$|I|=1$}
\begin{center}
  \includegraphics[scale=0.3]{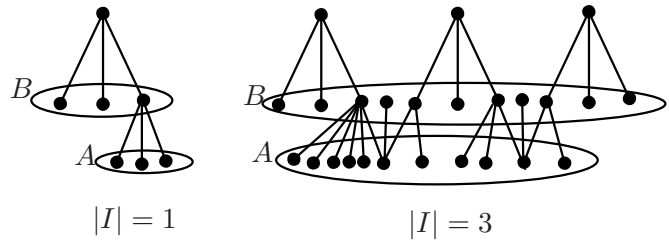}\\
\end{center}
\vspace{-3mm}
  \caption{The HIT $T_B$}\label{operation1}
\end{figure}

Since $(A,B)$ is $(\ve',\delta)$-super-regular,
for each $v\in \Gamma(x_i, B)$, $|\Gamma(v,A)|\ge \delta|A|$.
By condition (i), we have
$|\Gamma(x_i, B)|>\ve' |B|$ for each $i$.
Applying Fact~\ref{Typical},
we then know that there are at least $(\delta-\ve')|A|>|I|$ vertices from $\Gamma(v,A)$
typical to  $\Gamma(x_{i+1}, B)$ for each $1\le i\le |I|-1$.
That is, there exists $A_1\subseteq \Gamma(v,A)$ with $|A_1|\ge (\delta-\ve')|A|>|I|$
such that for each $a_1\in A_1$,  $|\Gamma(a_1, \Gamma(x_{i+1}, B))|\ge (\delta-\ve')d'|B|\ge 3|I|$.
%
As $ deg(x,B)\ge d'|B|\ge 3|I|$ for any $x\in I$ and
$(\delta-\ve')d'|B|\ge 3|I|$, combining the above argument,
we know
there is a claw-matching $M_I$ from $I$ to $B$ centered in $I$
such that  one vertex from $\Gamma(x_i, V(M_I))$
and one vertex from $\Gamma(x_{i+1}, V(M_I))$  have at least $(\delta-\ve')|A|>|I|$ common neighbors
in $A$.
Let $x_{i1}, x_{i2}, x_{i3}$ be the three neighbors of $x_i$ in $M_I$\,(in fact in $B$)
and suppose that $|\Gamma(x_{i3}, A)\cap \Gamma(x_{i+1,1}, A)|\ge |I|$.
For $1\le i\le |I|-1$,
we then choose distinct vertices $y_i\in \Gamma(x_{i3}, A)\cap \Gamma(x_{i+1,1}, A)$.
By  condition (iv),
there is a $\wedge$-matching $M_2$ between the vertex set $\{x_{i3}\,:\, 1\le i\le |I|-1 \}$
and the vertex set $A-\{y_i\,:\,1\le i\le |I|-1\}$ centered in the first set,
a matching  $M_3$ between $\{x_{i+1,1}\,:\, 1\le i\le |I|-1 \}$ and
$A-\{y_i\,:\,1\le i\le |I|-1\}-V(M_2)$ covering the first set, and a
matching $M_4$ between the vertex set $\{y_i\,:\,1\le i\le |I|-1\}$
and $B-V(M_I)$ covering the first set. Finally, by using (iv) again,
we can find three distinct
vertices $y_{31}, y_{32}, y_{33}\in \Gamma(x_{13}, A)-\{y_i\,:\,1\le i\le |I|-1\}-V(M_2)-V(M_3)$.
Let $T_{B}$ be the graph with
$$
V(T_{B})=V(M_I)\cup \{y_i\,:\,1\le i\le |I|-1\}\cup V(M_2)\cap V(M_3)\cup V(M_4)\cup \{y_{31}, y_{32}, y_{33}\}
$$
and
$$
E(T_{B})=M_I\cup \{y_ix_{i3}, y_ix_{i+1,1}\,:\,1\le i\le |I|-1\}\cup M_2\cup M_3\cup M_4\cup  \{x_{13}y_{31}, x_{13}y_{32}, x_{13}y_{33}\}.
$$
If $|I|=1$,  we choose $x_{11},x_{12}, x_{13}\in \Gamma(x_1,B)$ and $y_{31}, y_{32}, y_{33}\in \Gamma(x_{13},A)$.
Then let $T_B$ be the graph with
$$
  V(T_B)= \{x_1,x_{11}, x_{12}, x_{13}, y_{31}, y_{32}, y_{33}\}
  $$
 and
 $$
  E(T_B) = \{x_1x_{11}, x_{1}x_{12}, x_1x_{13}, x_{13}y_{31}, x_{13}y_{32}, x_{13}y_{33}\}.
$$
In any case, we see that $T_B$ is a HIT satisfying
\begin{eqnarray}\label{op1}
  |V(T_B)\cap B| &=& |V(T_B)\cap A| = 4|I|-1, \nonumber \\
  |L(T_B)\cap B| &=& \min\{2|I|+1, 3|I|-1\},\, |L(T_B)\cap A|=3|I|.
\end{eqnarray}
We call $T_B$ the insertion HIT associated with $B$.
Figure~\ref{operation1} gives a depiction of $T_B$ for $|I|=1,3$,  respectively.

\noindent \textbf{Operation II}\quad
Let $(A,B)$ be an $(\ve',\delta)$-super-regular pair, and $I$ a set of vertices disjoint from $A\cup B$. Suppose
that  (i) $ deg(x,A)\ge d'|A|> \ve'|A|$ and $ deg(x,A)\ge d'|A|\ge 3|I|$ for any $x\in I$;
(ii) $(\delta-\ve')d'|A|\ge 3|I|$; (iii) $(\delta-2\ve')|B|>|I|$; and
(iv) $\delta|B|>3|I|$.
Then we can do the following operations on $(A,B)$ and $I$.

Let $I=\{x_1,x_2,\cdots, x_{|I|}\}$.
We first assume that $|I|\ge 3$.

Since $(A,B)$ is $(\ve',\delta)$-super-regular,
for each $v\in \Gamma(x_i, A)$, $|\Gamma(v,B)|\ge \delta|B|$.
By condition (i), we have
$|\Gamma(x_i, A)|>\ve' |A|$ for each $i$.
Applying Fact~\ref{Typical},
we then know that there are at least $(\delta-2\ve')|B|>|I|$ vertices from $\Gamma(v,B)$
typical to  both $\Gamma(x_{i+1}, A)$ and $\Gamma(x_{i+2}, A)$  for each $1\le i\le |I|-2$.
That is, there exists $B_1\subseteq \Gamma(v,B)$ with $|B_1|\ge (\delta-2\ve')|B|>|I|$
such that for each $b_1\in B_1$,  $|\Gamma(b_1, \Gamma(x_{i+1}, A))|, |\Gamma(b_1, \Gamma(x_{i+2}, A))| \ge (\delta-\ve')d'|A|\ge 3|I|$.
%
As $ deg(x,A)\ge d'|A|\ge 3|I|$ for any $x\in I$ and
$(\delta-\ve')d'|A|\ge 3|I|$, combining the above argument,
we know
there is a claw-matching $M_I$ from $I$ to $A$ centered in $I$
such that  any one vertex from $\Gamma(x_i, V(M_I))$,
 any one vertex from $\Gamma(x_{i+1}, V(M_I))$, and any one vertex
 from $\Gamma(x_{i+2}, V(M_I))$  have at least $|I|$ common neighbors
in $B$.
Let $x_{i1}, x_{i2}, x_{i3}$ be the three neighbors of $x_i$ in $M_I$\,(in fact in $A$).
For $i=1$, choose $y_0\in \Gamma(x_{13},A)\cap \Gamma(x_{23},A)\cap \Gamma(x_{33},A)$.
Let $h=\lceil (|I|-3)/2\rceil$.
For $1\le k\le h$,
we then choose distinct vertices $y_k\in \Gamma(x_{1+2k,2}, A)\cap \Gamma(x_{2+2k,3}, A)\cap \Gamma(x_{3+2k,3}, A)$
\,(if $|I|=2+2k$, let $\Gamma(x_{3+2k,3}, A)=A$).
By  condition (iv),
there is a matching $M$ between the vertex set $\{x_{i3}, x_{1+2k,2}\,:\, 1\le i\le |I|, 1\le k\le h \}$
and the vertex set $B-\{y_0, y_k\,:\,1\le k\le h \}$ covering  the first set.
If $|I|$ is even, choose $y_{31}, y_{32}\in \Gamma(x_{13},B)$
such that they have not been chosen before; if
$|I|$ is odd, choose $y_{31}, y_{32}, y_{33}\in \Gamma(x_{13},B)$
such that they have not been chosen before.
%
Let $T_{A}$ be the graph with
$$
V(T_{A})=\left\{
  \begin{array}{ll}
    V(M_I)\cup V(M)\cup  \{y_0, y_k\,:\,1\le k\le h\}\cup  \{y_{31}, y_{32}\}, & \hbox{if $|I|$ is even;} \\
    V(M_I)\cup V(M)\cup  \{y_0, y_k\,:\,1\le k\le h\}\cup  \{y_{31}, y_{32}, y_{33}\}, & \hbox{if $|I|$ is odd;}
  \end{array}
\right.
$$
and $E(T_{A})$ containing all edges in $M_I\cup M\cup \{y_0x_{13},y_0x_{23},y_0x_{33}\} $ and all
edges in 
$$
\left\{
  \begin{array}{ll}
    \{x_{1+2k,2}y_k, x_{2+2k,2}y_k, x_{3+2k,2}y_k, x_{1+2h,2}y_h, x_{2+2h,2}y_h\,:\,1\le k\le h-1\}\cup  \{y_{31}, y_{32}\}, & \hbox{if $|I|$ is even;} \\
    \{ x_{1+2k,2}y_k, x_{2+2k,2}y_k, x_{3+2k,2}y_k\,:\,1\le k\le h\}\cup  \{y_{31}, y_{32}, y_{33}\}, & \hbox{if $|I|$ is odd.}
  \end{array}
\right.
$$
If $|I|=1$, we choose $x_{11},x_{12}, x_{13}\in \Gamma(x_1,A)$ and $y_{31}, y_{32}\in \Gamma(x_{13},B)$,
and then let $T_{A}$ be the graph with
$$
  V(T_{B}) = \{x_1,x_{11}, x_{12}, x_{13}, y_{31}, y_{32}\}\,\, \mbox{and}\,\,
  E(T_{B}) = \{x_1x_{11}, x_{1}x_{12}, x_1x_{13}, x_{13}y_{31}, x_{13}y_{32}\}.
$$
If $|I|=2$, we choose $x_{11},x_{12}, x_{13}\in \Gamma(x_1,A)$,
$x_{11},x_{12}, x_{13}\in \Gamma(x_2,A)$, $y\in \Gamma(x_{13},B)\cap \Gamma(x_{21},B)$,
$y_{11},y_{12}\in \Gamma(x_{13},B)$, and $y_{21},y_{22}\in \Gamma(x_{21},B)$
such that they are all distinct,
then let $T_{A}$ be the graph with
$$
  V(T_{B}) = \{x_i,x_{i1}, x_{i2}, x_{i3}, y, y_{i1}, y_{i2}\,:\, i=1,2\}\quad \mbox{and}\quad
$$
$$
  E(T_{B}) = \{x_ix_{i1}, x_ix_{i2}, x_ix_{i3}, x_{13}y, x_{21}y, x_{13}y_{11}, x_{13}y_{12}, x_{21}y_{21}, x_{21}y_{22}\}.
$$
We see that $T_A$ is a tree which has a degree 2 vertex $y$ only if $|I|=2$ and a degree
2 vertex $y_h$ only if $|I|>2$ and $|I|$ is even. In addition, $T_A$
satisfies the following.
\begin{eqnarray}\label{op2}
  |V(T_A)\cap A| &=&3|I|\quad \mbox{and}\quad   |L(T_A)\cap A|=\left\{
  \begin{array}{ll}
    2|I|, & \hbox{if $|I|= 1,2$;}  \nonumber \\
   2|I|-\left\lceil\frac{|I|-3}{2}\right\rceil, & \hbox{if $|I| \ge 3$;} \quad \mbox{and}\quad
  \end{array}
\right. \\
  |V(T_A)\cap B|&=&\left\{
  \begin{array}{ll}
    2, & \hbox{if $|I|= 1$;}  \nonumber \\
    2|I|+1, & \hbox{if $|I| \ge 2$;} \quad \mbox{and}\quad
  \end{array}
\right. \\
|L(T_A)\cap B|&=&\left\{
  \begin{array}{ll}
    2|I|, & \hbox{if $|I|= 1,2$;}  \\
    2|I|-\left\lceil\frac{|I|-3}{2}\right\rceil , & \hbox{if $|I| \ge 3$.}
  \end{array}
\right.
\end{eqnarray}
In this case, we call $T_A$ the insertion tree associated with $A$.
Notice that $|L(T_A)\cap A|=|L(T_A)\cap B|$ always holds.
Figure~\ref{operation1} gives a depiction of $T_A$ for $|I|=1,2,5,6$,  respectively.
\begin{center}
\begin{figure}[!htb]
\psfrag{B}{$B$} \psfrag{A}{$A$}
\psfrag{|I|=5}{$|I|=5$} \psfrag{|I|=6}{$|I|=6$}
\psfrag{|I|=1}{$|I|=1$}
\psfrag{|I|=2}{$|I|=2$}
\begin{center}
  \includegraphics[scale=0.3]{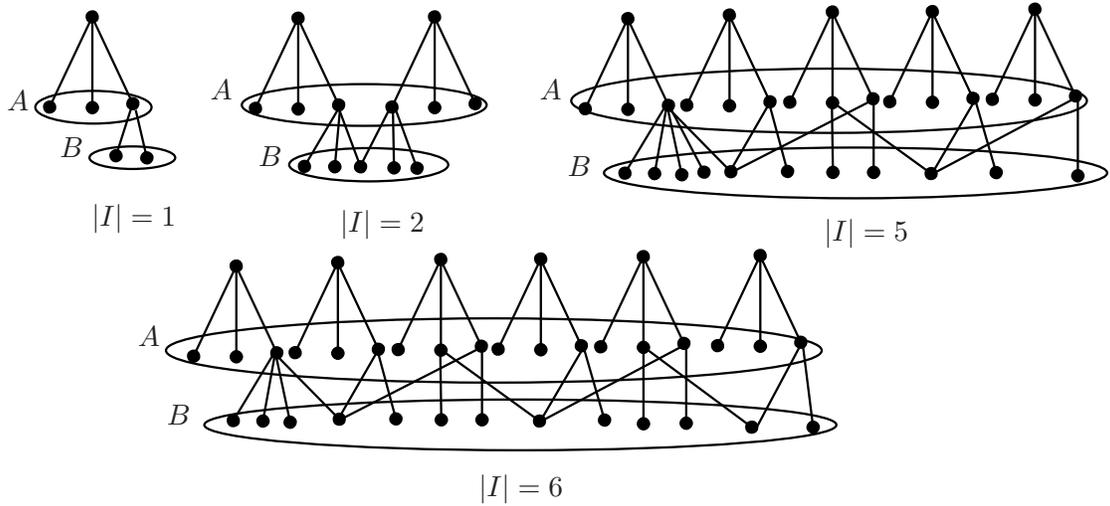}\\
\end{center}
\vspace{-2mm}
  \caption{The tree $T_A$}\label{operation2}
\end{figure}
\end{center}


\noindent \textbf{Operation III}\quad
Let $(B,F)$ be an $(\ve',\delta)$-super-regular pair, and $I$ a set of vertices disjoint from $B\cup F$. Suppose
that  $ deg(x,F)\ge d'|F|\ge 3|I|$ for any $x\in I$ and $\delta|B|\ge 6|I|$.
Then we can do the following operations on $(A,B)$ and $I$.

Let $I=\{x_1,x_2,\cdots, x_{|I|}\}$.
Since $ deg(x,B)\ge d'|B|\ge 3|I|$ for any $x\in I$,
there is a claw-matching $M_I$ from $I$ to $F$ centered in $I$.
Then as $\delta|B|\ge 6|I|$, there is a $\wedge$-matching
$M_{\wedge}$ from $V(M_I)\cap F$ to $B$ centered
in $V(M_I)\cap F$.
Let $T_{F}$ be the graph with
$$
V(T_{B})=V(M_I)\cup V(M_{\wedge})\quad \mbox{and} \quad
E(T_{B})=M_I\cup M_{\wedge}.
$$
We see that $T_{F}$ is a forest with no vertex of degree 2 satisfying
\begin{eqnarray}\label{op3}
  |V(T_F)\cap F|=|S(T_F)\cap F |&=& 3|I|\,\, \mbox{and}\,\, |V(T_F)\cap B| = |L(T_F)\cap B|=6|I|.
 \end{eqnarray}
 We call $T_F$ the insertion forest associated with $F$.

Now for each H-pair
$(X_i', Y_i')$, we may assume that $I(X_i')\ne \emptyset$
and $I(Y_i')\ne \emptyset$ for a uniform discussion,
as the consequent argument is independent of
the assumptions.
Recall that  $(X_i', Y_i')$ is $(2\ve, d-2\ve)$-super-regular by Fact~\ref{super-regularized XiYi}.
Notice that $deg(x, X_i') \ge (d-\ve)|X_i'|$ for each $x\in I(X_i')$,
$|I(X_i')|\le \frac{d^2N}{20}$, and $|X_i'|, |Y_i'|\ge (1-\ve)N$.
By simple calculations, we see that
(i) $deg(x,X_i')\ge (d-\ve)|X_i'|>2\ve|X_i'|$ and $(d-\ve)|X_i'|\ge 3d^2N/20 $ for each $x\in I(X_i')$; (ii)
$(d-2\ve-2\ve)(d-\ve)|X_i'|>3d^2N/20$; (iii) $(d-4\ve)|Y_i'|>d^2N/20$; and
(iv) $(d-2\ve)|Y_i'|>d^2N/5\ge 4I(X_i')$. Thus
all the conditions in
Operation I are satisfied.
So  we can find
a HIT $T_{X_i'}$ associated with $X_i'$.
As $|V(T_{X_i'})\cap X_i'|=|V(T_{X_i'})\cap Y_i'| \le 4|I(X_i')|\le \frac{d^2N}{5}$,
we know that $(X_i'-V(T_{X_i'}), Y_i'-V(T_{X_i'}))$ is $(4\ve, d-2\ve-d^2N/5)$-super regular.
Since $deg(y, Y_i')\ge  (d-\ve)|Y_i'|$ for each $y\in I(Y_i')$, we get
$deg(y, Y_i'-V(T_{X_i'}))\ge  (d-\ve-d^2/5)|Y_i'|$ for each $y\in I(Y_i')$.
By direct checking, conditions $(i)\sim (iv)$ of Operation I are satisfied
by the pair $(X_i'-V(T_{X_i'}), Y_i'-V(T_{X_i'}))$
and $I(Y_i')$.
Then we use Operation I  on $(X_i'-V(T_{X_i'}), Y_i'-V(T_{X_i'}))$
and $I(Y_i')$ to get a  HIT $T_{Y_i'}$ associated with $Y_i'-V(T_{X_i'})$.  Denote
$$
X_i^*=X_i'-V(T_{X_i'})-V(T_{Y_i'})  \quad \mbox{and} \quad Y_i^*=Y_i'-V(T_{X_i'})-V(T_{Y_i'}).
$$
By using~(\ref{op1})
in Operation I, we have $|X_i^*|=|Y_i^*|\ge (1-2d^2/5-\ve)N\ge N/2$. By Slicing lemma\,(Lemma~\ref{slicing lemma})
and Fact~\ref{super-regularized XiYi}, we have the following.
\begin{FAC}\label{ready-pair}
For each H-pair $(X_i,Y_i)$, $(X_i^*, Y_i^*)$ is  $(4\ve, d-2\ve-2d^2/5)$-super-regular with
$|X_i^*|=|Y_i^*|$. We call $(X_i^*, Y_i^*)$ a {\it ready H-pair}.
\end{FAC}

Then for each H-triple $(X_i', Y_i', F')$,
we may assume that $I(X_i')\ne \emptyset$
and $I(F')\ne \emptyset$\,(recall that $Y_i$ is assumed to be the dominator of
F, so $I(Y_i')=\emptyset$ by the distribution principle of vertices in $V_0'$ from Claim~\ref{absorbing}).
By Fact~\ref{super-regularized H-triple}, we know that
$(X_i', Y_i')$ is
$(2\ve, d-3.1\ve)$-
super-regular
and $(Y_i',F')$ is $(4.2\ve, d-3.1\ve-2d^3)$-super-regular.
Notice also that
$|X_i'|\ge (1-\ve)N$, $|Y_i'|\ge (1-3.1\ve)N$, $|F'|\ge (1-2.1\ve-2d^3)N/2$,
and $deg(x,X_i')\ge (d-\ve)|X_i'|$ and $deg(y,F')\ge (d-\ve)|F'|$ for each $x\in I(X_i')$
and each $y\in I(F')$.
Since $|I(X_i')|, |I(F')|\le \frac{d^2N}{20}$ and $\ve \ll d\ll 1$,
the conditions of Operation III are satisfied by $(Y_i', F')$ and $I(F')$
by direct calculations. Let $T_{F'}$ be the insertion forest associated with $F'$.
Then we
use Operation II  on $(X_i', Y_i'-V(T_{F'}))$ and $I(X_i')$ to get a  tree $T_{X_i'}$
associated with $X_i'$.  Denote
$$
X_i^*=X_i'-V(T_{X_i'}),\, Y_i^*=Y_i'-V(T_{F'})-V(T_{X_i'}),   \quad \mbox{and} \quad F^*=F'-V(T_{F'}).
$$
By using~(\ref{op2}) and~(\ref{op3}) in Operation II and Operation III, respectively,
we have $|X_i^*|, |Y_i^*|\ge (1-3.1\ve-9d^2/20)N\ge N/2$
and $|F^*|\ge  (1-2.1\ve-2d^3)N/2-3d^2N/20\ge (1-2.1\ve-2d^3-3d^2/10)N/2$.
By Slicing lemma
and Fact~\ref{super-regularized H-triple}, we have the following.

\begin{FAC}\label{ready-triple}
For each H-triple $(X_i,Y_i,F)$,
$(X_i^*,Y_i^*)$ is  $(4\ve, d-3.1\ve-9d^2/20)$-super-regular and
$(Y_i^*,F^*)$ is $(8.4\ve, d-2.1\ve-3d^2/10-2d^3)$-super-regular.
We call $(X_i^*, Y_i^*, F^*)$ a {\it ready H-triple}.
\end{FAC}

\noindent \textbf{Step 5.} Apply the Blow-up lemma to find a HIT within each ready H-pair and among each ready H-triple.

In order to apply the Blow-up Lemma,
we first give two lemmas which assure the existence of a given subgraph in a complete bipartite graph.

\begin{LEM}\label{complte-bipartite}
Suppose $0<\ve\ll d\ll 1$ and $N$ is a large integer.
If $G(A,B)$ is a balanced complete bipartite graph with $(1-\ve-d^2/2)N\le |A|=|B|\le N$,
then $G(A,B)$ contains a HIST $T_{pair}$ with $\Delta(T_{pair})\le \lceil 2/d^3\rceil$
and $||L(T_{pair})\cap A|-|L(T_{pair})\cap B||=\ell $ for any given non-negative integer $\ell$ with  $ \ell \le d^2N$.
\end{LEM}

\pf  By the symmetry, we only show that we can construct a HIST $T$ such that $|L(T)\cap A|-|L(T)\cap B|=\ell$.
Let $\Delta'=\lceil d^3N\rceil$.
We choose distinct $a_1,a_2,\cdots, a_{\Delta'}\in A$ and distinct $b_1,b_2,\cdots, b_{\Delta'-1}\in B$.
Then we decompose all vertices in $B$ into $B_1, B_2,\cdots, B_{\Delta'}$ such that  $3\le |B_i|\le 1/d^3$,
$B_i\cap B_{i+1}=\{b_i\}$ for  $1\le i \le \Delta'-1$, and $B_i\cap B_j =\emptyset$ for $|i-j|>1$.
Now we choose $\ell+1$  distinct vertices $b_{\Delta'}, b_{\Delta'+1},\cdots, b_{\Delta'+\ell}$
from $B-\{b_i\,:\,1\le i \le \Delta'-1 \}$.
As  $\Delta'=\lceil d^3N\rceil$, $\ell+\Delta'\le (d^2 +d^3)N+1$, and thus
\begin{equation*}
   2(\ell+\Delta') \le  (2d^2 +2d^3)N+2 \le (1-d^2/2-\ve)N- \lceil d^3N\rceil\le |A|-\lceil d^3N \rceil.
\end{equation*}
Thus we can use all of the vertices in $\{b_i\,:\,1\le i\le \Delta'+\ell\}$ to cover all
vertices in $A-\{a_i\,|\,1\le i\le \Delta'-1\}$ such that each $b_i$ can be
adjacent to at least two distinct vertices. We partition $A-\{a_i\,|\,1\le i\le \Delta'-1\}$
arbitrarily into $A_1,A_2, \cdots, A_{\ell+\Delta'}$ such that $2\le |A_i|\le 1/d^3 $.
Now let $T$ be a spanning subgraph of $G(A,B)$ such that
$$
E(T)=\{a_ib \,|\, b\in B_i, 1\le i\le \Delta'  \}\cup \{b_ja\,|\, a\in A_j, 1\le j\le \Delta'+\ell\}.
$$
Clearly, $\Delta(T)\le \lceil 2/d^3\rceil$. As $|A|=|B|$,
$|S(T)\cap A|=\Delta'$,  and $|S(T)\cap B|=\Delta'+\ell$, we then have that
$|L(T)\cap A|-|L(T)\cap B|=\ell$.  We denote $T$ as $T_{pair}$.
 \qqed

\begin{LEM}\label{half-complte-tripartite}
Suppose $0<\ve\ll d\ll 1$ and $N$ is a large integer.
Let $G=G(A,B,F)$ be a tripartite graph with $V(G)$ partitioned into $A\cup B\cup F$  such that
  both $G[A\cup B]$  and
$G[B\cup F]$ are complete bipartite graphs. If
(i) $(1-4\ve-d^2/2)N  \le |A|, |B|\le N$, (ii) $(1/2-2.1\ve-3d^2/20-d^3)N\le |F|\le (1/2-d^3)N$, and (iii)
for any given non-negative integer $l\le 3d^2N/10$, we have
$|B|-2(|A\cup F|-|B|-l)\ge 3d^3N/2$ holds, then
 $G$ contains a HIST $T_{triple}$ and a path $P_{triple}$ spanning on a subset of $L(T_{triple})$ such that
\vspace{-4mm}
\begin{itemize}
  \item [$($a$)$] $T_{triple}$ is a HIST of $G$ with $\Delta(T_{triple})\le \lceil3/d^3\rceil$;
  \item [$($b$)$] $|L(T_{triple})\cap B|=|L(T_{triple})\cap (A\cup F)|-l$.
  \item [$($c$)$]  $P_{triple}$ is a $(b,f)$-path on $L(T_{triple})\cap F$ and  any $|L(T_{triple})\cap F|$ vertices from $L(T_{triple})\cap B$, and $|V(P_{triple})\cap F|\le 5d^2N/6$.
 \end{itemize}
\vspace{-4mm}
\end{LEM}

\pf Let $\Delta'=\lceil d^3N/2\rceil$.
We choose distinct $b_1,b_2,\cdots, b_{\Delta'}\in B$  and
partition all vertices in $F$ into $F_1, F_2,\cdots, F_{\Delta'}$ such that $3\le |F_i|\le 1/d^3$.
Then we choose  distinct $a_1,a_2,\cdots, a_{\Delta'-1}\in A$  and
decompose all vertices in $A$ into $A_1, A_2,\cdots, A_{\Delta'}$ such that $3\le |A_i|\le 2/d^3$,
$A_i\cap A_{i+1}=\{a_i\}$ for $1\le i \le \Delta'-1$,  and  $A_i\cap A_j =\emptyset$ for $|i-j|>1$.
Choose one more vertex, say $a_{\Delta'}\in A-\{a_i\,|\, 1\le i \le \Delta'-1\}$.
Let $l'=|A\cup F|-|B|-l$. Notice that $l'>0$.
Now we choose $l'$  distinct vertices $f_1, f_2,\cdots, f_{l'}$ from $A-\{a_i\,:\, 1\le  i\le \Delta'\}\cup F$\,(choose as many as possible from $F$ first)
and partition any $2l'$ vertices of $B-\{b_i\,:\, 1\le i\le \Delta'\}$ into $B_1, B_2,\cdots, B_{l'}$
such that $|B_i|=2$. By (iii), we see that there are at least
$\lfloor d^3N\rfloor$ vertices left in
$B'=B-\{b_i\,:\, 1\le i\le \Delta'\}-\bigcup_{i=1}^{l'}\{B_i\}$.
Hence we can partition  $B'=B_1'\cup B_2'\cup \cdots \cup B'_{\Delta'}$ such that $|B'_{\Delta'}|\ge 2$ and $|B_j'|\ge 1$
for $j\ne \Delta'$.  We let $T$ be a subgraph of $G$ on $A\cup B\cup F$ with
$$
E(T)=\{b_if, b_ia,a_ib' \,:\, f\in F_i,a\in A_i,b'\in B'_i, 1\le i\le \Delta' \}\cup \{f_ib\,:\, b\in B_i, 1\le i\le l'\}.
$$
By the construction, $T$ is a HIST of $G$, which clearly satisfies (a).
Since $|S(T)\cap B|=\Delta'$ and $|S(T)\cap (A\cup F)|=\Delta'+l'=\Delta'+|A\cup F|-|B|-l$,
we then see that $T$ satisfies (b).
If  $L(T)\cap F\ne \emptyset$, let $f\in L(T)\cap F$ and $b\in L(T)\cap B$, we can then take a $(b,f)$-path $P$
with $V(P)\cap F=L(T)\cup F$ and $|V(P)|=2|L(T)\cap F|$.  By (i) and (ii),
we see that $l'=|A\cup F|-|B|-l\ge (1/2-6.1\ve-4d^2/5-d^3)N$.
Hence
$|V(P)\cap F|=|F|-l'\le 5d^2N/6$.
Denote $T$ as $T_{triple}$ and $P$ as $P_{triple}$.
\qqed

Now for $1\le i\le t$ and for each ready H-pair $(X_i^*,Y_i^*)$, suppose, without of loss generality, that
$|(L(T_{X_i'})\cap Y_i')\cup ((L(T_{Y_i'})\cap Y_i') |-
|((L(T_{X_i'})\cap X_i')\cup ((L(T_{Y_i'})\cap X_i') |=l'$,
where  $T_{X_i'}$ is the insertion HIT associated with $X_i'$ and $T_{Y_i'}$ is the insertion HIT
associated with $Y_i'$. Notice that  $l'\le d^2N$ from~(\ref{op1}) and~(\ref{op2}).
Let $x_a\in S(T_{X_i'})\cap X_i'$ be a non-leaf of $T_{X_i'}$ and
$y_b\in S(T_{Y_i'})\cap Y_i'$  a non-leaf of $T_{Y_i'}$.
Since $(X_i', Y_i')$ is $(2\ve, d-2\ve)$-super-regular by Fact~\ref{super-regularized XiYi}
and $|Y_i'-Y_i^*|\le 2d^2N/5$, we have
$deg(x_a, Y_i^*)\ge (d-2\ve-d^2/2)N\ge dN/2 $. Similarly,
$deg(y_b, X_i^*)\ge (d-2\ve-d^2/2)N\ge dN/2 $.
Also, from Step 2,  we have $\Gamma(x_i^*, Y_i), \Gamma(x_i^{**}, Y_i)\ge (d-3\ve)N$.
So, $\Gamma(x_i^*, Y^*_i), \Gamma(x_i^{**}, Y^*_i)\ge (d-3\ve-d^2/2)N \ge dN/2$.
Similarly, we have
$\Gamma(y_i^*, X^*_i), \Gamma(y_i^{**}, X^*_i)\ge (d-3\ve-d^2/2)N \ge dN/2$.
Recall that $(X_i^*,Y_i^*)$ is $(4\ve, d-2\ve-8d^2/20)$-super-regular by Fact~\ref{ready-pair},
and therefore $(X_i^*,Y_i^*)$ is $(4\ve, d/2)$-super-regular.
By the the strengthened version of the Blow-up lemma
and Lemma~\ref{complte-bipartite}\,(the conditions
are certainly satisfied by $X_i^*$ and $Y_i^*$), we can find a HIST $T_{1}^i\cong T_{pair}$ on $X_i^*\cup Y_i^*$ such that there exist
$y_a\in S(T^i_1)\cap \Gamma(x_a, Y_i^*)$, $x_b\in S(T^i_1)\cap \Gamma(y_b, X_i^*)$,
$y_i'\in S(T^i_1)\cap \Gamma(x_i^*, Y_i)$, $y_i''\in S(T^i_1)\cap \Gamma(x_i^{**}, Y_i)$,
and $x_i'\in S(T^i_1)\cap \Gamma(y_i^*, X_i)$, $x_i''\in S(T^i_1)\cap \Gamma(y_i^{**}, X_i)$
such that $|L(T_1^i)\cap X_i^*|-|L(T_1^i)\cap Y_i^*|=l'$.
Hence
$|L(T_1^i)\cap X_i^*|+|L(T_{X_i'})\cap X_i'|+|L(T_{Y_i'})\cap X_i'|=|L(T_1^i)\cap Y_i^*|+|L(T_{X_i'})\cap Y_i'|+|L(T_{Y_i'})\cap Y_i'|$.
Let $T^i=T^i_1\cup T_{X_i'}\cup T_{Y_i'}\cup\{x_ay_a, y_bx_b\}\cup \{x_i^*y_i',x_i^{**}y_i'', y_i^*x_i', y_i^{**}x_i''\}$. It is clear that
$T^i$ is a HIST on $X'_i\cup Y_i'\cup I(X_i')\cup I(Y_i')$ such that
$$
\{x_i^*,x_i^{**},y_i^*, y_i^{**}\}\subseteq L(T^i)\quad \mbox{and} \quad |L(T^i)\cap X_i'|=|L(T^i)\cap Y_i'|.
$$
For the ready H-pair $(X_0^*,Y_0^*)$,
let $x_a\in S(T_{X_0'})\cap X_0'$ be a non-leaf of $T_{X_0'}$ and
$y_b\in S(T_{Y_0'})\cap Y_0'$  a non-leaf of $T_{Y_0'}$.
By the the strengthened version of the Blow-up lemma
and Lemma~\ref{complte-bipartite}\,(the conditions
are certainly satisfied by $X_0^*$ and $Y_0^*$), we can find a HIST $T_{1}^0\cong T_{pair}$ on $X_0^*\cup Y_0^*$ such that there exist
$y_0'\in S(T^0_1)\cap \Gamma(x_0^*, Y_0)$, $y_0''\in S(T^0_1)\cap \Gamma(x_0^{**}, Y_0)$,
$x_{t+1}'\in S(T^0_1)\cap \Gamma(y_{t+1}^*, Y_0)$, $x_{t+1}''\in S(T^0_1)\cap \Gamma(y_{t+1}^{**}, Y_0)$,
and $x_0'\in S(T^0_1)\cap \Gamma(y_0^*, X_0)$, $x_0''\in S(T^i_1)\cap \Gamma(y_0^{**}, X_0)$
such that
$|L(T_1^0)\cap X_0^*|+|L(T_{X_0'})\cap X_0'|+|L(T_{Y_0'})\cap X_0'|=|L(T_1^0)\cap Y_0^*|+|L(T_{X_0'})\cap Y_0'|+|L(T_{Y_0'})\cap Y_0'|+2$.
Let $T^0=T^0_1\cup T_{X_0'}\cup T_{Y_0'}\cup\{x_ay_a, y_bx_b\}\cup \{x_0^*y_0',x_0^{**}y_0'', y_0^*x_0', y_0^{**}x_0'',
y_{t+1}^*x_{t+1}', y_{t+1}^{**}x_{t+1}''\}$. It is clear that
$T^0$ is a HIST on $X'_0\cup Y_0'\cup I(X_0')\cup I(Y_0')$ such that
$$
\{x_0^*,x_0^{**},y_0^*, y_0^{**}, y_{t+1}^*, y_{t+1}^{**}\}\subseteq L(T^0)\quad \mbox{and} \quad |L(T^0)\cap X_0'|=|L(T^i)\cap Y_0'|+2.
$$


For each ready triple  $(X_i^*,Y_i^*, F^*)$,
we know that $(X_i^*,Y_i^*)$ is  $(4\ve, d-3.1\ve-9d^2/20)$-super-regular
and $(Y_i^*,F^*)$ is $(8.4\ve, d-2.1\ve-3d^2/10-2d^3)$-super-regular
by Fact~\ref{ready-triple}.
Notice that $(1-4\ve-9d^2/20)N\le |X_i^*|, |Y_i^*|\le N$ and
$(1/2-2.1\ve-3d^2/30-d^3)N\le |F^*|\le (1/2-d^3)N$.
Let $|I(X_i')|=l'$ and $|I(F')|=l/6$ for some integer $l$.
By Operation II we have $|V(T_{X_i'})\cap X_i'|\le 3l'$ and $|V(T_{X_i'})\cap Y_i'|\le 2l'+1$.
By Operation III we have
$|V(T_{F'})\cap F_i'|=l/2$ and $|V(T_{F'})\cap Y_i'|=l$.
Notice that $|L(T_{X_i'})\cap X_i'|=|l(T_{X_i'})\cap Y_i'|$.
Hence,
\begin{eqnarray*}
  |Y_i^*|-2(|X_i^*\cup F^*|-|Y_i^*|-l) &\ge &3(|Y_i'|-2l'-l-1)-2(|X_i'|-3l')-2(|F'|-l/2)+2l \\
   &=& 3|Y_i'|-2|X_i'|-2|F'|-3 \\
   &\ge & 3(1-3.1\ve)N-2N-N+2d^3N-3>3d^3N/2.
\end{eqnarray*}
By the weak version of the Blow-up lemma\,(Lemma~\ref{blow-up})
and Lemma~\ref{half-complte-tripartite}, we then can find a HIT $T_1^i\cong T_{triple}$ on $X_i^*\cup Y_i^*\cup F^*$
and a path $P_i\cong P_{triple}$ spanning on $L(T_1^i)\cap F^*$ and other $|L(T_1^i)\cap F^*|$ vertices from $Y_i^*$.
Let $y_a\in S(T_{X_i'})\cap Y_i'$ be a non-leaf of $T_{X_i'}$\,(take $y_a$ as the degree 2 vertex if $T_{X_i'}$ has one) and
 $y'_a\in S(T_{F'})\cap Y_i'$  a non-leaf of $T_{F'}$.
Then as
$(Y_i',F')$ is $(4.1\ve, d-2.1\ve-2d^3)$-super-regular,
we have $|\Gamma(y_a, F')|, |\Gamma(y'_a, F')| \ge (d-2.1\ve-2d^3)N/2 $.
Since  $|F'-F^*|\le 3d^2N/20$, we then know that
$|\Gamma(y_a, F^*)|, |\Gamma(y'_a, F^*)| \ge (d-2.1\ve-3d^2/10-2d^3)N/2$.
Since $|F^*\cap L(T_1^i)|=|V(P_i)\cap F^*|\le 5d^2N/6<(d-2.1\ve-3d^2/10-2d^3)N/2$,
there exist $f_a\in (S(T_1^i)\cap F^*)\cap \Gamma(y_a, F^*)$ and $f'_a\in (S(T_1^i)\cap F^*)\cap \Gamma(y'_a, F^*)$.
For each $x\in I(F')$, since $deg(x, F')\ge (d-\ve)|F'|\ge (d-\ve)(1-2.1\ve-d^3)N/2$,
we know there exists $f'\in (S(T_1^i)\cap F^*)\cap \Gamma(x,F^*)$.
From Step 2, we have $|\Gamma(x_i^*, Y_i)\cap \Gamma(x_i^{**}, Y_i)|\ge (d-\ve)^2N$ and
$|\Gamma(y_i^*, X_i)\cap \Gamma(y_i^{**}, X_i)|\ge (d-\ve)^2N$.
Hence $|\Gamma(x_i^*, Y'_i)\cap \Gamma(x_i^{**}, Y'_i)|\ge ((d-\ve)^2-3.1\ve)N$.
Since $|S(T_1^i\cup T_{X_i'}\cup T_{F'})\cap X_i'|<d^2N/2$,
we see that there exists $y'\in \Gamma(x_i^*, Y_i)\cap \Gamma(x_i^{**}, Y_i)\cap L(T_1^i\cup T_{X_i'}\cup T_{F'})$.
Similarly, there exists $x'\in \Gamma(y_i^*, X_i)\cap \Gamma(y_i^{**}, X_i)\cap L(T_1^i\cup T_{X_i'}\cup T_{F'})$.
Let $T^i=T^i_1\cup T_{X_i'}\cup T_{F'}\cup\{xf'\,:\, x\in I(F'), f'\in (S(T_1^i)\cap F^*)\cap \Gamma(x,F^*) \}\cup\{y_af_a, y'_af'_a\}\cup \{y'x_i^*, y'x_i^{**}, x'y_i^*, x'y_i^{**}\}$.
It is clear that
$T^i$ is a HIST on $X'_i\cup Y_i'\cup F'\cup I(X_i')\cup I(F')$
such that
$$
\{x_i^*,x_i^{**},y_i^*, y_i^{**}\}\subseteq L(T^i)\quad \mbox{and}\quad |L(T^i)\cap X_i'|=|L(T^i)\cap Y_i'|.
$$
Let $H^i=T^i\cup P_i$.
We call $P_i$ the accompany path of $T^i$.

\noindent \textbf{Step 6.} Apply the Blow-up Lemma again on the regular-pairs induced on the leaves of each HIT obtained in Step 5
to find two vertex-disjoint paths covering all the leaves. Then connect all the HITs into a HIST of $G$ and connect the
disjoint paths into a cycle  using the edges initiated in Step 2.

Suppose $1\le i\le t$.
For each H-pair $(X_i,Y_i)$, let $X_i^L=X_i'\cap L(T^i)-\{x_i^*,x_i^{**}\}$ and $Y_i^L=Y_i'\cap L(T^i)-\{y_i^*,y_i^{**}\}$,
and for each H-triple $(X_i,Y_i,F)$, let  $X_i^L=X_i'\cap L(T^i\cup P_i)-\{x_i^*,x_i^{**}\}$
and $Y_i^L=Y_i'\cap L(T^i\cup P_i)-\{y_i^*,y_i^{**}\}$,
where $T^i$ is the HIST found in Step 5, and $P_i$ is the
accompany path of $T^i$.
By Operations I, II and III, and
the proofs of the Lemmas~\ref{complte-bipartite} and~\ref{half-complte-tripartite},
we have
$I(X_i')\cup I(Y_i')\subseteq S(T_i)$
and $F'\cup I(F')\subseteq S(T_i\cup P_i)$. Thus,
 $X_i^L\cup Y_i^L=L(T^i)-\{x_i^*,x_i^{**}, y_i^{*}, y_i^{**}\}$ for each H-pair
 and  $X_i^L\cup Y_i^L=L(T^i\cup P_i)-\{x_i^*,x_i^{**}, y_i^{*}, y_i^{**}\}$
for each H-triple.  Furthermore, we have $|X_i^L|=|Y_i^L|$.
For the H-pair $(X_0,Y_0)$, let $X_0^L=X_0'\cap L(T^0)-\{x_0^*,x_0^{**}, y_{t+1}^{*}, y_{t+1}^{**}\}$ and $Y_0^L=Y_0'\cap L(T^0)-\{y_0^*, y_0^{**}\}$.
We have $X_0^L\cup Y_0^L=L(T^0)-\{x_0^*,x_0^{**}, y_{t+1}^{*}, y_{t+1}^{**}\}$ and $|X_0^L|=|Y_0^L|$
since from Step 5 we have $|L(T^0)\cap X_0'|=|L(T^0)\cap Y_0'|+2$.
By the construction of
$T_{pair}$ and $H_{triple}$, we see that $|S(T_i)\cap X_i'|, |S(T_i)\cap Y_i'|\le d^2N$.
Since each  H-pair $(X_i',Y_i')$ is $(2\ve, d-2\ve)$-super-regular, and
each pair $(X_i',Y_i')$ from an H-triple $(X_i',Y_i',F')$ is
$(2\ve, d-3.1\ve)$-
super-regular, by Slicing Lemma, we then know that
$(X_i^L, Y_i^L)$ is $(4\ve, d-3.1\ve-d^2)$-super-regular and hence is
$(4\ve, d/2)$-super-regular.

For each $1\le i\le t$, by the choice of
$x_i^*,x_i^{**}, y_i^*, y_i^{**}$, we have
$|\Gamma(x_i^*,Y_i)|, |\Gamma(x_i^{**},Y_i)|\ge (d-\ve)N$ and $|\Gamma(y_i^*,X_i)|,|\Gamma(y_i^{**},X_i)|\ge (d-\ve)N$.
Hence, $|\Gamma(x_i^*,Y_i^L)|, |\Gamma(x_i^{**},Y_i^L)|\ge (d-\ve-3.1\ve-d^2)N>dN/2$ and $|\Gamma(y_i^*,X_i^L)|,|\Gamma(y_i^{**},X_i^L)|\ge (d-\ve-3.1\ve-d^2)N>dN/2$.
Similar results hold for the vertices $x_0^*,x_0^{**}, y_{t+1}^*, y_{t+1}^{**}$.
For each $0\le i\le t$, we choose distinct vertices $y_i'\in \Gamma(x_i^*,Y_i^L)$, $y_i''\in \Gamma(x_i^{**},Y_i^L)$ and
$x_i'\in \Gamma(y_i^*,X_i^L)$, $x_i''\in \Gamma(y_i^{**},X_i^L)$.
If $T^i$ does not have the accompany path,
then by the strengthened version of the Blow-up lemma, we can find an
$(x'_i, y_i')$-path $P_1^i$ and an $(x''_i, y_i'')$-path $P_2^i$ such that $P_1^i\cup P_2^i$ is spanning on $X^L_i\cup Y^L_i $.
If $T^i$ has the accompany $(b,f)$-path $P_i$, we see that $deg(b, X_i^L), deg(f, Y_i^L) \ge dN/2$
as $(X_i', Y_i')$ is
$(2\ve, d-3.1\ve)$-
super-regular,
and $(Y_i',F')$ is $(4.2\ve, d-3.1\ve-2d^3)$-super-regular.
Applying  the strengthened version of the Blow-up lemma, we can find an
$(x'_i, y_i')$-path $P_{11}^i$ and an $(x''_i, y_i'')$-path $P_2^i$ such that $P_{11}^i\cup P_2^i$ is spanning on $X^L_i\cup Y^L_i $,
and two consecutive internal vertices $a', b'\in V(P_{11}^i)$ with $b'\in \Gamma(f, Y_i^L) $, and $a'\in \Gamma(b, X_i^L)$.
Let $P^i_{1}=P_{11}^i\cup P_i\cup \{fb', ba'\}-\{a'b'\}$. Notice that
for the H-pair $(X_0, Y_0)$, the two vertices $y_{t+1}^*, y_{t+1}^{**}$ are
not used in this step, but we will connect them to $y_0^*$ and $y_0^{**}$,
respectively, in next step.

We now connect the small HITs and paths together to
find an SGHG of $G$.
In Case A, for $1\le i\le t-1$, we have $|S(T^i)\cap Y_i|\ge d^3N/2>\ve N$ and $|S(T^{i+1})\cap X_{i+1}|\ge d^3N/2>\ve N$. Since
$(Y_i, X_{i+1})$ is an $\ve$-regular pair with density $d$, we see that there is an edge $e_i$ connecting
$S(T^{i+1})\cap X_{i+1}$ and $S(T^{i+1})\cap X_{i+1}$. Let
$$
T=\bigcup_{i=1}^{t}T^i\cup \{e_i\,|\, 1\le i\le t-1\}.
$$
Then $T$ is a HIST of $G$.
Let $C$ be the cycle formed by all the paths in $\bigcup_{i=1}^{t}(P_1^i\cup P_2^i)$ and
all edges in the following set
$$
\{x_i^*y_i', x_i^{**}y_i^{''},
y_i^*x_{i}',y_i^{**}x_{i}^{''}\,: 1\le i\le t\}\cup \{y_i^*x_{i+1}^*,y_i^{**}x_{i+1}^{**}\,:1\le i\le t-1\}\cup \{y_t^*x_1^{**}, y_t^{**}x_1^*\},
$$
notices that the edges in $\{y_i^*x_{i+1}^*,y_i^{**}x_{i+1}^{**}\,:1\le i\le t-1\}\cup \{y_t^*x_1^{**}, y_t^{**}x_1^*\}$
above are guaranteed in Step 2.
It is easy to see that $C$ is a cycle on $L(T)$.
Hence $H=T\cup C$ is an SGHG of $G$.

In Case B, for $1\le i\le t-1$, we have $|S(T^i)\cap Y_i|\ge d^3N/2>\ve N$ and $|S(T^{i+1})\cap X_{i+1}|\ge d^3N/2>\ve N$. Since
$(Y_i, X_{i+1})$ is an $\ve$-regular pair with density $d$, we see that there is an edge $e_i$ connecting
$S(T^{i+1})\cap X_{i+1}$ and $S(T^{i+1})\cap X_{i+1}$.  Similarly, there is an edge $e_0$ connecting $S(T_0)\cap X_0$ and
$S(T^{1})\cap X_{1}$.
Let
$$
T=\bigcup_{i=1}^{t}T^i\cup \{e_i\,|\, 0\le i\le t-1\}.
$$
Then $T$ is a HIST of $G$.
Let $C$ be the cycle formed by all paths in $\bigcup_{i=1}^{t}(P_1^i\cup P_2^i)$ and
all edges in the set  $\{y_0^*y_{t+1}^*, y_0^{**}y_{t+1}^{**}, y_{t+1}^*x_1^*, y_{t+1}^{**}x_1^{**}, x_0^*y_t^{**}, x_0^{**}y_t^*\}$ and
in the  following set
$$
\{x_i^*y_i', x_i^{**}y_i^{''},
y_i^*x_{i}',y_i^{**}x_{i}^{''}\,: 0\le i\le t\}\cup\{y_i^*x_{i+1}^*, y_i^{**}x_{i+1}^{**}\,:1\le i\le t-1\}.
$$
%
%
It is easy to see that $C$ is a cycle on $L(T)$.
Hence $H=T\cup C$ is an SGHG of $G$.

The proof of Theorem~\ref{non-extremal} is now finished.
\qqed

\subsection{Proof of Theorem~\ref{extremal_1}}

%

By the assumption that $deg(v_1,V_2)\leq 2\beta n$  for each $v_1\in V_1$
and the assumption that $\delta(G)\geq(2n+3)/5$ in Extremal Case 1, we see that
\begin{equation}\label{mindregreeG[v1]}
 \delta(G[V_1])\geq (2n+3)/5-2\beta n.
\end{equation}
Then (\ref{mindregreeG[v1]}) implies that
\begin{equation}\label{V1V2size}
  |V_1|\geq (2n+3)/5-2\beta n \quad \mbox{and} \quad |V_2|\leq 3n/5+2\beta n.
\end{equation}

Also, by $|V_2|\geq (2/5-4\beta)n$ in the assumption,
\begin{equation}\label{V1size}
  |V_1|\leq (3/5+4\beta)n.
\end{equation}

We will construct an SGHG of $G$ following several steps below.

\noindent \textbf{Step 1. Repartitioning}

Set  $\alpha_1=\alpha^{1/3}$ and $\alpha_2=\alpha^{2/3}$.
Let
$$
V_1'=V_1\quad \mbox{and}\quad V_2'=\{v\in V_2\,|\, deg(v,V_1)\leq \alpha_1|V_1|\}.
$$
 Then by $d(V_1,V_2)\leq \alpha $, we have
 $$
\alpha_1|V_1||V_2-V_2'|\le e(V_1,V_2')+e(V_1,V_2-V_2')=e(V_1,V_2)\leq\alpha|V_1||V_2|.
$$
This gives that
\begin{equation}\label{V2-V2 prime size}
    |V_2-V_2'|\leq \alpha_2|V_2|.
\end{equation}
Denote $V_{12}^0=V_2-V_2'$. Then by the definition of $V_2'$, we have
\begin{equation}\label{mindegreev120v1}
  \delta(V_{12}^0,V_1')>\alpha_1|V_1'| \quad \mbox{and}\quad  \delta(G[V_2'])\geq (2n+3)/5-\alpha_1|V_1'|\geq(2/5-\alpha_1(3/5+4\beta))n,
\end{equation}
where the last inequality follows from (\ref{V1size}).

Let $n_i=|V_i'|$ for $i=1, 2$. Then by (\ref{mindregreeG[v1]}) and (\ref{V1size}),
\begin{equation}\label{mindregreeGv1prime 2}
     \delta(G[V_1'])\geq (2n+3)/5-2\beta n \geq \frac{2/5-2\beta }{3/5+4\beta}n_1\geq (2/3-8\beta)n_1,
\end{equation}
and by (\ref{V1V2size}) and the second inequality in (\ref{mindegreev120v1}),
\begin{eqnarray}\label{mindregreeG[v2]prime}
                  \delta(G[V_2'])&\geq &(2/5-\alpha_1(3/5+4\beta))n \ge \frac{(2/5-\alpha_1(3/5+4\beta))}{3/5+2\beta}n_2 \ge (2/3-1.1\alpha_1)n_2, \nonumber
\end{eqnarray}
provided that $\beta\leq\frac{0.3\alpha_1}{9\alpha_1+20/3}$.

\noindent \textbf{Step 2. Finding three connecting edges}

AS $G$ is 3-connected, there are 3 independent edges $x_L^1y_L^1,x_L^2y_L^2$
and $x_Ny_N$ connecting $V_1'\cup V_{12}^0$ and $V_2'$
such that $x_L^1,x_L^2,x_N\in V_1'\cup V_{12}^0$ and $y_L^1,y_L^2 , y_N\in V_2'$.
In the remaining steps, we will find a HIST $T_1$ in $G[V_1'\cup V_{12}^0]$
 with $x_N$ as a non-leaf and $x_L^1,x_L^2$ as leaves,
and a HIST $T_2$ of $G[V_2']$
 with $y_N$ as a non-leaf and $y_L^1,y_L^2$ as leaves.
Then $T=T_1\cup T_2\cup\{x_Ny_N\}$ is a HIST of $G$.
 By finding a hamiltonian $(x_L^1,x_L^2)$-path $P_1$ on $L(T_1)$,
and a hamiltonian $(y_L^1,y_L^2)$-path on $L(T_2)$, we see that
$$
C:=P_1\cup P_2\cup \{x_L^1y_L^1, x_L^2y_L^2\}
$$
forms a cycle on $L(T)$. Hence $H:=T\cup C$ is an SGHG of $G$.

\noindent \textbf{Step 3. Initiating two HITs}

In this step, we first initiate a HIT in  $G[V_1'\cup V_{12}^0]$
containing $X_N$ as a non-leaf and $x_L^1$ and $x_L^2$ as leaves.
 Then, we initiate a HIT in $G[V_2']$ containing $y_N$ as a non-leaf
and $y_L^1$ and $y_L^2$ as leaves.

For $x_L^1,x_L^2,x_N\in V_1'\cup V_{12}^0$, by (\ref{mindregreeG[v1]}) and (\ref{mindegreev120v1}),
each of them has at least $\alpha_1|V_1'|\geq 9$ neighbors in $V_1'$.
Thus, we choose distinct $z_L^1,z^1,z_L^2,z^2,z_N^1,z_N^2,z_N^3 \in V_1'$ such that
$$
    x_L^1\sim z_L^1,z^1, \quad x_L^2\sim z_L^2,z^2, \quad x_N\sim z_N^1,z_N^2,z_N^3.
$$
(Note that $x_L^1$ and $x_L^2$ may be from $V_{12}^0$, and therefore they may not have too
many neighbors in $V_1'$, we then choose $z_L^1$ and $z_L^2$ from $V_1'$ as their neighbors, respectively.)

\begin{center}
\begin{figure}[!htb]
\psfrag{B}{$B$} \psfrag{A}{$A$}
\psfrag{z_1^1}{$z_1^1$} \psfrag{x_L^1}{$x_L^1$}
\psfrag{z_L^1}{$z_L^1$}
\psfrag{z^1}{$z^1$} \psfrag{z^11}{$z^{11}$}
\psfrag{z_2^2}{$z_2^2$} \psfrag{x_L^2}{$x_L^2$}
\psfrag{z_L^2}{$z_L^2$}
\psfrag{z^2}{$z^2$} \psfrag{z^22}{$z^{22}$}
\psfrag{z^12}{$z^{12}$} \psfrag{v_1^R}{$v_1^R$}
\psfrag{x_N}{$x_N$} \psfrag{z_N^1}{$z_N^1$}
\psfrag{z_N^2}{$z_N^2$}
\psfrag{z_N^3}{$z_N^3$}
\psfrag{z_N^11}{$z_N^{11}$}
\begin{center}
  \includegraphics[scale=0.3]{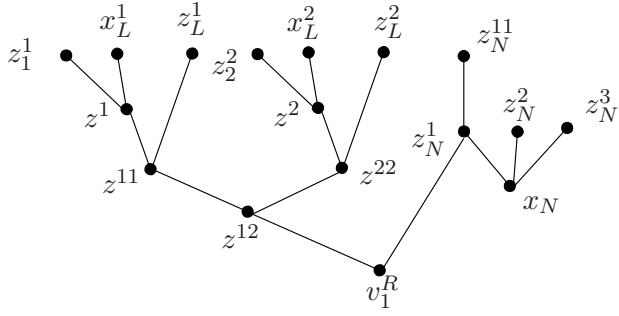}\\
\end{center}
\vspace{-2mm}
  \caption{The tree $T_{11}$}\label{t11}
\end{figure}
\end{center}

By (\ref{mindregreeGv1prime 2}),
we see that any two vertices in $G[V_1']$ have at least $(1/3-16\beta) n_1\geq 14$ neighbors in common.
Thus, we can choose distinct vertices $z^{11},z^{22},z^{12},v_1^{R} \in V_1'-\{x_L^1,x_L^2,x_N,z_L^1,z^1,z_L^2,z_N^1,z_N^2,z_N^3\}$
such that
$$
    z^{11}\sim z_L^1,z^1, \quad z^{22}\sim z_L^2,z^2, \quad z^{12}\sim z^{11},z^{22},\quad v_1^{R}\sim z^{12},z_N^1.
$$
Furthermore, by (\ref{mindregreeGv1prime 2}) again,
we have $\delta(G[V_1'])\geq (2/3-8\beta)n_1\geq 17$.
Choose $z_1^1, z_2^2, z_N^{11}\in V_1'$ not chosen above such that
$$
z_1^1\sim z^1, z_2^2\sim z^2, z_N^{11}\sim z_N^1.
$$
Let $T_{11}$ be the graph with
$$
V(T_{11})=\{x_L^1,x_L^2,x_N,z_N^1,z_L^1,z^1,z_L^2,z^{11},z^{12},z^{22},z^{2},z_N^2,z_N^3,v_1^{R},z_1^1, z_2^2, z_N^{11}\}
$$
and with edges indicated above except the edges $x_L^1z_L^1$ and $x_L^2z_L^2$.
We see that $T_{11}$ is a tree with $v_1^{R}$ as the only degree 2 vertex, and
$|V(T_{11})|=17$ and $ |L(T_{11})|=9$. Notice that
in $T_{11}$, $z_L^1,x_L^1$ and $z_L^2,x_L^2$ are leaves, and $x_N$ is a non-leaf.
Figure~\ref{t11} gives a depiction of $T_{11}$.

Notice that the edges $x_L^1z_L^1$ and $x_L^2z_L^2$ are not used in $T_{11}$.
We will first construct a HIST $T_1$ in  $G[V_1^1\cup V_{12}^0]$
containing  $T_{11}$  as a subgraph,  then find a
hamiltonian $(z_L^1,z_L^2)$-path on $L(T_1)-\{x_L^1,x_L^2\}$ by Lemma~\ref{Hamiltonian-connected},
finally by adding $x_L^1z_L^1$ and $x_L^2z_L^2$ to the path,
we get a hamiltonian  $(x_L^1,x_L^2)$-path on $L(T_1)$.
The reason that we avoid using $x_L^1$ and $x_L^2$ is that when $x_L^1,x_L^2\in  V_{12}^0 $,
we may not be able to have the condition of Lemma~\ref{Hamiltonian-connected} on $G[L(T_1)]$
in our final construction.

Then we initiate a HIT in $G[V_2']$ containing $y_L^1, y_L^2$ as leaves, and $y_N$ as a non-leaf.

As $y_L^1, y_L^2, y_N\in V_2'$, by (\ref{mindregreeG[v2]prime}) and the fact that
each two vertices from $V_2'$ have at least $(1/3-2.2\alpha_1)n_2\geq 7$
common neighbors implied from (\ref{mindregreeG[v2]prime}), we can choose distinct vertices
$$
y^{12},y_N^1,y_N^2,y_N^3, v_2^{R}\in V_2'-\{y_L^1,y_L^2,y_N\}
$$
such that
\begin{equation}\label{y12}
y^{12}\sim y_L^1, y_L^2, \quad y_N\sim y_N^1,y_N^2,y_N^3,\quad  v_2^{R}\sim y^{12},y_N.
\end{equation}
Let $T_{21}$ be  the graph with $V(T_{21})=\{y_L^1,y_L^2,y_N, y^{12}, y_N^1,y_N^2,y_N^3,v_2^{R}\}$
and with $E(T_{21})$ described as in (\ref{y12}).

We see that $T_{21}$ is a tree with $v_2^{R}$ the only degree 2 vertex
and $y_L^1,y_L^2 \in L(T_{21})$, $y_N\in S(T_{21})$ and
\begin{equation}\label{VTILTI size}
|V(T_{21})\cap V_2'|=8, \quad |L(T_{21})\cap V_2'|=5.
\end{equation}
Denote
$$
U_1=V_1'-V(T_{11}), \quad U_2=V_2'-V(T_{21}),\quad \mbox{ and} \quad V_{12}=V_{12}^0-V(T_{11}).
$$

\noindent \textbf{Step 4. Absorbing vertices in $V_{12}^0$}

We may assume that $V_{12}^0\neq\emptyset$. For otherwise, we skip this step.
Let $|V_{12}|=n_{12}$ and  $V_{12}^0=\{x_1,x_2, \cdots, x_{n_{12}}\}$.

Since $|V(T_{11})|=17,$ by (\ref{mindegreev120v1}), we get
$$
\delta(V_{12}^0,U_1)> \alpha_1|V_1'|-17\geq 3\alpha_2|V_2|\geq 3|V_2-V_2'|\geq 3|V_{12}^0|.
$$
Thus, there is a claw-matching $M_c$ from $V_{12}^0$ to $U_1$ centered in $V_{12}^0$.
For $i=1,2,\cdots,n_{12}$, let $x_{i1},x_{i2}$ and $x_{i3}$  be the three neighbors of $x_i$ in $M_c$.
If $n_{12}=1$,
let $T_a=M_c$, and we finish this step.
Thus we assume $n_{12}\geq 2$.

By (\ref{mindregreeGv1prime 2}), each two vertices in in $V_1'$ have at least
\begin{equation}\label{V1primeneighbor}
    (1/3-16\beta)n_1\geq 6\alpha_2|V_{12}^0|+17
\end{equation}
neighbors in common. The above inequality holds as $n_1\ge 2n/5-2\beta n$,
$|V_2|\le 3n/5+2\beta n$ by (\ref{V1V2size}), and we can assume that $18\alpha_2/5+106\beta/15+12\alpha_2\beta+18/n-32\beta^2\le 2/15$.

Thus, for each $i=1,2,\cdots,n_{12}-1$, we can find  distinct vertices $x_{13}^i,x_{23}^i,x_{i3}^3,x_{i+1,1}^3$ in $U_1-V(M_c)$
such that
\begin{equation}\label{edges-absorbing}
  x_{13}^i\sim x_{i3}, x_{i+1,1}, \quad x_{23}^i\sim x_{13}^i,\quad x_{i3}^3\sim x_{i3}, \quad x_{i+1,1}^3 \sim x_{i+1,1}.
\end{equation}
Let $T_a$ be the graph with
$V(T_a)=V(M_c)\cup \{x_{13}^i,x_{23}^i,x_{i3}^3,x_{i+1,1}^3\,:\,1\leq i\leq n_{12}-1\}$,
and $E(T_a)$ including all edges indicated in (\ref{edges-absorbing}) for all $i$ and all edges in $M_c$.
It is easy to see, by the construction,  that $T_a$ is a HIT with
$$
|V(T_a)\cap U_1|=7n_{12}-4 \quad \mbox{and}\quad |L(T_a)\cap U_1|=4n_{12}-1.
$$

Using (\ref{V1primeneighbor}) again, we can find $x_N^{11}\in U_1-V(T_a)$ such that
$x_N^{11}\sim v_1^R, x_{11}$, where $v_1^R\in V(T_{11})$ and $ x_{11}\in V(T_a)$. By (\ref{mindregreeGv1prime 2}),
$$
\delta[G[V_1']]\geq (2n+3)/5-2\beta n\geq 6\alpha_2|V_{12}^0|+20,
$$
since $|V_2|\leq 2n/5+2\beta n$, and we can assume that
$2\beta-12\alpha_2\beta-18\alpha_2/5-21/n\le 2/5$.
So we can find distinct vertices $x_N^{12}, x_{11}^1\in U_1-V(T_a)-\{x_N^{11}\}$
such that $x_N^{12}\sim x_N^{11}, x_{11}^{1}\sim x_{11}$.

Let $T_1^1$ be the graph with
$$
V(T_1^1)=V(T_{11})\cup V(T_a) \cup \{ x_N^{11},x_N^{12}, x_{11}^{1}\} \quad \mbox{and}\quad
E(T_1^1)=E(T_{11})\cup E(T_a) \cup \{ x_N^{11}v_R^1,x_N^{11}x_{11}, x_N^{12}x_N^{11}, x_{11}^{1}x_{11}\}.
$$
Then  $T_1^1$ is a HIT such that
\begin{equation}\label{T12size}
|V(T_1^1)\cap U_1|=7n_{12}+16 \quad \mbox{and}\quad |S(T_1^1)\cap U_1|=3n_{12}+7.
\end{equation}
Denote
$U_1'=U_1-V(T_1^2)$ and $ U_2'=U_2-V(T_1^2)$.

\noindent \textbf{Step 5. Completion of HITs $T_1$ and $T_2$}

In this step, we construct a HIST $T_i$ in $G[V_i']\, (i=1,2)$ containing $T_1^i$ as an induced subgraph.

The following lemma guarantees the existence of a specified HIST in a graph with $n$
vertices and minimum degree at least $(2/3-\alpha')n $ for some $0<\alpha'\ll 1.$

\begin{LEM}\label{TH construction}
Let $H$ be an $n$-vertex graph with $\delta(H)\geq(2/3-\alpha')n$ for some constant $0<\alpha'\ll 1.$ Then $H$ has a HIST $T_H$ satisfying
\vspace{-4mm}
\begin{itemize}
  \item [$($i$)$]$T_H$ has a vertex $v_R$ of degree at least $(2/3-\alpha')n-1$, and $v_R$ can be chosen arbitrarily from $V(H)$;
  \item [$($ii$)$]$|S(T_H)|\le (1/6+\alpha'/2)n+2$.
\end{itemize}
\end{LEM}

\pf Let $v_R\in V(H)$ be an arbitrary vertex. If $n(mod\,2)\equiv deg(v_R)+1(mod\,2)$, then we let $N_R=N_H(v_R)$.
For otherwise, let $N_R$  be a subset of $N(v_R)$ with $|N_H(v_R)|-1$ elements.
Let $T_{v_R}$ be the star with $V(T_{v_R})=\{{v_R}\}\cup N_R $ and $E(T_R)=E(\{v_R\},N_R)$.
Let $V_0=V(H)-V(T_{v_R})$. By $\delta(H)\geq(2/3-\alpha')n$, we have $|V_0|\leq( 1/3+\alpha')n+1$.
By the choice of $N_R$, we have $|V_0|\equiv0(mod \,2)$.
If $V_0=\emptyset$, then let $T_H=T_{v_R}$. For otherwise, we claim as follows.

\begin{CLA}\label{residues-absorbing}
Let $V_1\subseteq V(H)$ be a subset with $|V_1|\geq(2/3-\alpha')n-1$ and $|V_1|(mod \,2)\equiv n(mod \,2)$.
Then there exist two vertices from $V_0=V(H)-V_1$ such that they have a common neighbor in $V_1$.
\end{CLA}

\proof[Proof of Claim~\ref{residues-absorbing}]
 We assume that $|V_1|\le(2/3+2\alpha')n$. For otherwise, $|V_0|<(1/3-2\alpha')n$.
Since $\delta(H)\geq (2/3-\alpha')n$,
any two vertices of $H$ have at least $(1/3-2\alpha')n$ neighbors in common.
By $|V_0|<(1/3-2\alpha')n$, any two vertices from $V_0$ have a common neighbor from $V_1$. We are done.
Thus $|V_1|\le(2/3+2\alpha')n$, and hence $|V_0|\ge (1/3-2\alpha')n\geq 3$.
By the assumption that $|V_1|\geq (2/3-\alpha')n-1$, we have $|V_0|\le (1/3+\alpha')n+1$.
This implies that $\deg(v_0,V_1)\geq(1/3-2\alpha')n-2$ for each $v_0\in V_0$.
As $|V_0|\geq 3$ and $3((1/3-2\alpha')n-2)>(2/3+2\alpha')n>\ge |V_1|$\,(provided that $8\alpha'+6/n< 1/3$),
we see that there must be two vertices from $V_0$ such that they have a neighbor in common in $V_1$.
\qed

By Claim~\ref{residues-absorbing},
there exist two vertices $v_0^{11},v_0^{12}\in V_0$ such that they have a common neighbor in $T_{v_R}$.
Adding $v_0^{11}$ and $v_0^{12}$ to $T_{v_R}$ and two
edges connecting them to one of their common neighbor in $V(T_{v_R})$.
Let $T^1_{v_R}$ be the resulting graph. Then we see that
$T^1_{v_R}$ is a HIT with $|V(T^1_{v_R})|=|V(T_{v_R})|+2$, and hence $(|V(T_{v_R})|+2)(mod \,2)\equiv n(mod \,2)$.
Also  $|V(T^1_{v_R})|\geq |V(T_{v_R})|\geq (2/3-\alpha')n-1$.
So we can use Claim~\ref{residues-absorbing} again to find another pair of vertices from $V_0-\{v_0^{11}, v_0^{12}\}$
such that they have a common neighbor in $V(T^1_{v_R})$.
Adding the new pair of vertices and two edges connecting them to one of their common neighbor
in $V(T^1_{v_R})$ into $T^1_{v_R}$,
we get a new HIT $T^2_{v_R}$.
By repeating the above process another $l_0=(|V_0|-4)/2$ times,
we get a HIT $T^{l_0}_{v_R}$.
Let $T_H=T^{l_0}_{v_R}$. We claim that $T_H$ has the required properties in Lemma~\ref{TH construction}.
Notice first that $d_{T_H}(v_R)\geq (2/3-\alpha')n-1$.
Then since $T_H$ has $v_R$ and at most $|V_0|/2$ distinct vertices as non-leaves
and $|V_0|\leq (1/3+\alpha')n+1$, we see that $|S(T_H)|\le (1/6+\alpha'/2)n+2.$
\qed

Let $H_1=G[U_1'\cup \{v^1_R\}]$. Recall that $v^1_R$ is a non-leaf in $T_1^1$.
By (\ref{mindregreeGv1prime 2}) and (\ref{T12size}), and
by noticing that $n_{12}\le |V_2-V_2'|\le \alpha_2|V_2|\le 3\alpha_2n_1/2 $
(by (\ref{V1V2size})), we see that
\begin{eqnarray}
  \delta(H_1)  &\geq&(2/3-8\beta)n_1-(7n_{12}+19)  \nonumber\\
    &\geq& (2/3-8\beta)n_1-21\alpha_2n_1/2-19 \nonumber\\
     &\geq& (2/3-11\alpha_2)|V(H_1)|.
\end{eqnarray}

Let $\alpha'=11\alpha_2\ll 1$\, (by assuming $\alpha\ll 1$).
By Lemma~\ref{TH construction},
we can find a HIT $T_1'$ in $H_1$ with $v_R^1$ as the prescribed vertex in condition(i).
It is easy to see that $T_1:=T_1^1\cup T_1'$ is a HIST of $G[V_1'\cup V_{12}^0]$ and
\begin{eqnarray}\label{Ln2 size}
   s_1&=&|S(T_1)\cap V_1'|= |S(T^1_1)\cap V_1'|+|S(T_1')\cap V_1'|\nonumber \\
      &\leq & 3n_{12}+7+(1/6+5.5\alpha_2)|V(H_1)|+2 \,(\mbox{by}\,(\ref{T12size})\,\mbox{and Lemma}~\ref{TH construction}) \nonumber \\
  &\leq&  3n_{12}+9+(1/6+5.5\alpha_2)n_1\nonumber \\
      &\leq& (1/6+10.5\alpha_2)n_1 \,(\mbox{by}\, n_{12}\le 3\alpha_2n_1/2).
\end{eqnarray}

Let $H_2=G[U_2'\cup \{v_{R}^2\}]$. By(\ref{mindregreeG[v2]prime}) and (\ref{VTILTI size}),
we see that
\begin{eqnarray*}
  \delta(H_2)&\geq & (2/3-1.1\alpha_1)n_2-8\ge (2/3-1.2\alpha_1) |V(H_2)|.
 \end{eqnarray*}
By letting $\alpha'=1.2\alpha_1$,
we can find a HIT $T_2'$ in $H_2$ with $v_R^2$ as the prescribed
vertex in condition (i) of Lemma~\ref{TH construction}.
Then $T_2:=T_1^2\cup T_2' $ is a HIST of $G[V_2']$.
Also, notice that
\begin{eqnarray} \label{l^n2}
   s_2&=&|S(T_2)\cap V_2'|=|S(T^2_1)\cap V_2'|+|S(T_2')\cap V_2'|\nonumber \\
       &\leq & 3+(1/6+0.6\alpha_1)|V(H_2)|+2\nonumber \\
   &\leq& (1/6+0.7\alpha_2)n_2,
\end{eqnarray}
where the last inequality holds by assuming $5/n_2\leq 0.1\alpha_2$.

\noindent \textbf{Step 6. Finding two long paths.}

In this step, we first find a hamiltonian $(z_L^1,z_2^2)$-path $P_1^1$ in $G[L(T_1)-\{x_L^1,x_L^2\}]$;
then find a hamiltonian $(y_L^1,y_L^2)$-path $P_2$ in $G[L(T_2)]$.
Let $G_{11}=G[L(T_1)-\{x_L^1,x_L^2\}]$ and $n_{11}=|V(G_{11})|$.
We will show that $\delta(G_{11})>\frac{1}{2}n_{11}$.
We may assume $s_1\geq (1/6-8\beta)n_1-2$. For otherwise, if $s_1<(1/6-8\beta)n_1-2$,
 then by (\ref{mindregreeGv1prime 2}), we get
 \begin{eqnarray}
    \delta(G_{11})&\geq & \delta(G[V_1'])-s_1-2 \nonumber\\
    &\geq & (2/3-8\beta)n_1-((1/6-8\beta)n_1-1-2)-2 \nonumber\\
 &\geq & \frac{1}{2}n_1+1 \ge \frac{1}{2}n_{11}+1.\nonumber
 \end{eqnarray}
 Hence, $s_1\geq (1/6-8\beta)n_1-2$, implying that
 \begin{equation}\label{}
   n_{11}\leq (5/6+8\beta)n_1+2 \quad \mbox{and thus} \quad n_1\geq \frac{n_{11}-2}{5/6+8\beta}.
 \end{equation}
Hence, by (\ref{Ln2 size})
\begin{eqnarray}
\delta(G_{11})&\geq & \delta(G[V_1'])-s_1-2 \ge (2/3-8\beta)n_1-(1/6+10.5\alpha_2)n_1-2\nonumber\\
   &\geq & (1/2-8\beta-11\alpha_2)n_1\ge \frac{1/2-2\beta-11\alpha_2}{5/6+2\beta} (n_{11}-2)> n_{11}/2, \nonumber
  \end{eqnarray}
the last inequality holds by assuming $3\beta+11\alpha_2+2/n_{11}<1/12$.
By applying Lemma~\ref{TH construction} on $G_{11}$, we find a
hamiltonian $(z_L^1,z_L^2)$-path $P_1^1$ in $G_{11}$.
 Let $P_1=P_1^1\cup \{z_L^1x_L^1,z_L^2x_L^2\}$.
We see  that $P_1$ is a a hamiltonian $(x_L^1,x_L^2)$-path on $L(T_1)$.

Let $G_{22}=G[L(T_2)]$ and $n_{22}=|V(G_{22})|$.
We will show that $\delta(G_{22})> n_{22}/2$.
We may assume that  $s_2\geq (1/6-1.1\alpha_1)n_2-2$.
For otherwise, if $s_2<(1/6-1.1\alpha_1)n_2-2$,
then by (\ref{mindregreeG[v2]prime}), we see that
\begin{eqnarray*}
   \delta(G_{22})&\geq & \delta(G[V_2'])-s_2-2 \\
   &>& (2/3-1.1\alpha_1)n_2-((1/6-1.1\alpha_1)n_2-2)-2 \\
   &>& n_2/2\ge n_{22}/2.
\end{eqnarray*}

Thus, $s_2\geq (1/6-1.1\alpha_1)n_2-2$, implying that
$$
n_{22}\leq n_1-s_2\leq (5/6+1.1\alpha_1)n_2+2 \quad \mbox{ gives that} \quad n_2\geq \frac{n_{22}-2}{5/6+1.1\alpha_1}.
$$
By (\ref{mindregreeG[v2]prime}) and (\ref{l^n2}),
\begin{eqnarray}
  \delta(G_{22})&\geq & \delta(G[V_2'])-s_2-2 \nonumber  \\
  &\geq& (2/3-1.1\alpha_1)n_2-(1/6+0.7\alpha_1)n_2-2\nonumber  \\
  &\geq& (1/2-1.9\alpha_1)n_2\ge \frac{1/2-1.9\alpha_2}{5/6+1.1\alpha_2}(n_{22}-2) \nonumber \\
   &>& n_{22}/2.\nonumber
\end{eqnarray}
The last inequality follows by assuming that $2.45\alpha_1+2/n_{11}<1/12$.
Hence, by Lemma~\ref{TH construction},
there is a  hamiltonian $(y_L^1,y_L^2)$-path $P_2$ in $G_{22}$.

\noindent \textbf{Step 7. Forming an SGHG}

Let $T=T_1\cup T_2\cup \{x_Ny_N\}$ and $C=P_1\cup P_2\cup \{x_L^1y_L^1,x_L^2y_L^2\}$.
We see that $T$ is a HIST of $G$ with $L(T)=V(P_1)\cup V(P_2)$ and $C$ is a cycle spanning on $L(T)$.
Hence $H=T\cup C$ is an SGHG of $G$.


\subsection{Proof of Theorem~\ref{extremal_2}}

 Notice that the assumption of Extremal Case 2 implies that
 $$
|V_1|>(3/5-\alpha)n \quad \mbox{and} \quad    |V_2|\ge (2/5-2\beta)n.
 $$

 We may assume that the graph $G$ is minimal with respective to the number of edges. This implies
that no two adjacent vertices both have degree larger than $(2n+3)/5$.
 (For otherwise, we could delete any edges incident to two vertices both with degree larger than $(2n+3)/5$.)
 We construct an SGHG in $G$ step by step.

\noindent \textbf{Step 1. Repartitioning}

Set $\alpha_1=\alpha^{1/3} \quad \mbox{and} \quad \alpha_2=\alpha^{2/3}.$
Let
\begin{eqnarray*}
  V_2' &=& \{v\in V_2\,|\, deg(v, V_1)\geq (1-3\alpha_1)|V_1|\}, \\
  V_0' &=& \{v\in V_2-V_2'\,| \,deg(v, V_1)\leq\alpha_1|V_2|/6\}, \\
  V_1' &=& V_1\cup V_0', \quad  V_{12}^{0}\,=\,V_2-V_2'-V_0'.
\end{eqnarray*}
As $d(V_1, V_2)\geq 1-3\alpha$, the following holds,
\begin{eqnarray*}
     (1-3\alpha)|V_1||V_2|&\leq & e_{G}(V_1,V_2)=e_{G}(V_1,V_2')+e_{G}(V_1,V_2-V_2')\\
        &\le & |V_1||V_2'|+(1-3\alpha_1)|V_1||V_2-V_2'|.
\end{eqnarray*}

 The inequality implies that
 \begin{equation}\label{size of V2-v2'}
    |V_2-V_2'|\leq \alpha_2|V_2|.
 \end{equation}

 As a consequence of moving vertices in $V_2-V_2'$ out from $V_2$, by~(\ref{size of V2-v2'})
we get
 \begin{eqnarray}\label{degree-V1toV2'}
   \delta(V_1,V_2')&\geq & (2n+3)/5-2\beta n-\alpha_2|V_2|\nonumber\\
        &\geq & (2n+3)/5-6\beta |V_2|-\alpha_2|V_2|\nonumber\\
        &\geq & (2n+3)/5- 2\alpha_2|V_2|,
 \end{eqnarray}
provided that $6\beta \leq \alpha_2$.
And as a consequence of moving vertices in $V_0'$ to $V_1$,
 \begin{eqnarray}\label{degreeV0'toV2'}
   \delta(V_0', V_2') &\geq & \delta(G)-\Delta(V_0', V_1)-\Delta(V_0', V_2-V_2')\nonumber\\
    &\geq & (2n+3)/5-\alpha_1|V_2|/6-\alpha_2|V_2|\nonumber\\
    &\geq&(2n+3)/5-\alpha_1|V_2|/3 \,(\mbox{provided that\, $\alpha_2 \leq \alpha_1/6$}),
 \end{eqnarray}
 and
 \begin{equation}\label{V120toV1' mindeg}
    \alpha_1|V_2|/6<\delta(V_{12}^0, V_1')< (1-3\alpha_1)|V_1|.
\end{equation}
 From (\ref{degree-V1toV2'}) and (\ref{degreeV0'toV2'}), we have
 \begin{equation}\label{degree-V1'toV2'}
    \delta(V_{1}',V_2')\geq (2n+3)/5-\alpha_1|V_2|/3.
 \end{equation}
As
\begin{equation}\label{v2V1 mindeg}
    \delta(V_2',V_1') \geq (1-3\alpha_1)|V_1| \geq (1-3\alpha_1)(3/5-\alpha)n > \left\lceil(2n+3)/5 \right\rceil,
\end{equation}
we get that
\begin{equation}\label{degree-V1'}
    deg(v_1')=\left\lceil(2n+3)/5 \right\rceil
\end{equation}
for each $v_1' \in V_1'$, by the minimality assumption of $e(G)$.
Hence (\ref{degree-V1'toV2'}) and (\ref{degree-V1'}) give that
\begin{equation}\label{degree-maxV1'}
    \Delta(G[V_1'])\leq \alpha_1|V_2|/3.
\end{equation}

\noindent \textbf{Step 2. Finding a vertex $v_{\scr2}^{\scr *}$ from $V_2'$ with large degree in $V_1'$}

Let
\begin{equation}\label{ein}
    e_{in}=e(G[V_1'])
\end{equation}
be the number of edges within $V_1'$, notice that $e_{in}$ maybe 0.
Then
\begin{equation}\label{eg(v1v2v12)}
    e_G(V_1',V_2'\cup V_{12}^0) =|V_1'|\lceil(2n+3)/5\rceil -2 e_{in}.
\end{equation}
Let
\begin{equation}\label{din}
    d_{in}=e_{in}/|V_1'| \quad \mbox{and} \quad |n_0|=|\,|V_2'\cup V_{12}^0|-\lceil(2n+3)/5 \rceil|.
\end{equation}
By (\ref{degree-maxV1'}) and the definition of $d_{in}$ in (\ref{din}),
we have
\begin{equation*}\label{lower-din-larger}
    \lfloor d_{in}\rfloor\le \alpha_1|V_2|/6.
\end{equation*}
In fact, since $\Delta(V_1,V_1')\le \Delta(V_1,V_1)+\Delta(V_1,V_0')\le 2\beta n +|V_0'|\le 2\beta n+\alpha_2|V_2|$,
and $\Delta(V_0', V_1')\le \alpha_1|V_2|/6+\alpha_2|V_2|$, more precisely, we have
\begin{eqnarray}\label{lower-din-smaller}
   2d_{in}&= & 2e_{in}/|V_1'|\le (2\beta n+\alpha_2|V_2|)|V_1|/|V_1'|+(\alpha_1|V_2|/6+\alpha_2|V_2|)|V_0'|/|V_1'| \nonumber\\
     &\leq & (2\beta n+\alpha_2|V_2|)+\alpha_2(\alpha_1|V_2|/6+\alpha_2|V_2|)\,\,(\mbox{as $|V_0'|\le \alpha_2|V_2|$ and $|V_1|, |V_2|\le |V_1'|$}) \nonumber\\
 &\leq & (6\beta +\alpha_2+\alpha/6+\alpha_2^2)|V_2| \,\,(\mbox{as $\beta n\le 3\beta |V_2|$})\nonumber\\
  &\le & 2\alpha_2|V_2| \, (\mbox{provided that $6\beta+\alpha/6+\alpha_2^2\le \alpha_2$}).
\end{eqnarray}
Set

\noindent\textbf{Case A}.
    $\left\lceil(2n+3)/5  \right\rceil-|V_2'\cup V_{12}^0|=n_0\geq 0;$

\noindent \textbf{Case B}.
    $|V_2'\cup V_{12}^0|-\left\lceil(2n+3)/5 \right\rceil=n_0\geq 1.$

We have
\begin{numcases}{n_0=}
\lceil(2n+3)/5 \rceil- |V_2'\cup V_{12}^0|\leq 2\beta n+\alpha_2|V_2|\le  (6\beta+\alpha_2)|V_2|\leq 2\alpha_2|V_2|,& Case A, \nonumber\\
& \label{n0-upper-bound} \\
  |V_2'\cup V_{12}^0|-\lceil(2n+3)/5 \rceil \leq (2/5+\alpha)n -\lceil(2n+3)/5 \rceil\leq \alpha n,& Case B.\nonumber
\end{numcases}

Then in case A,
\begin{eqnarray*}
                   e_G(V_1',V_2'\cup V_{12}^0)  &=&|V_1'|\lceil(2n+3)/5 \rceil-2e_{in}\quad (\mbox{by~(\ref{degree-V1'})})\nonumber \\
                   &=& |V_1'|(|V_2'\cup V_{ 12}^0|+n_0-2d_{in})\nonumber \\
                   &\geq & |V_2'\cup V_{12}^0|(|V_1'|+1.4n_0-3.2d_{in} ),
\end{eqnarray*}
as
$1.4|V_2'\cup V_{12}^0|\leq 1.4((2n+3)/5+\alpha n) \leq (3/5-\alpha)n<|V_1'|$ and
$1.6|V_2'\cup V_{12}^0|\geq 1.6((2n+3)/5-2\beta-\alpha_2)n\geq (3/5+2\beta+\alpha_2)n)>|V_1'|$
provided that  $2.4\alpha <1/25$ and $5.2\beta+2.6\alpha_2\leq 1/25$ respectively.
Since $e_G(V_1',V_2'\cup V_{12}^0)\le |V_2'\cup V_{12}^0||V_1'|$, we have $|V_1'|+1.4n_0-3.2d_{in}\leq |V_1'|$,  and thus $1.4n_0\leq 3.2 d_{in}$.

In Case B,
 \begin{eqnarray*}
                   e_G(V_1',V_2'\cup V_{12}^0)  &=&|V_1'|\lceil(2n+3)/5 \rceil-2e_{in}\quad (\mbox{by~(\ref{degree-V1'})})\nonumber \\
                   &=& |V_1'|(|V_2'\cup V_{12}^0|-n_0-2d_{in})\nonumber \\
                   &\geq & |V_2'\cup V_{12}^0|(|V_1'|-1.6n_0-3.2d_{in} ),
               \end{eqnarray*}
as $
1.6|V_2'\cup V_{12}^0|\geq 1.6((2n+3)/5-2\beta-\alpha_2n) \geq (3/5+2\beta+\alpha_2)n>|V_1'|$
provided that  $5.2\beta+2.6\alpha_2\leq 1/25$.

Let
\begin{numcases}{d_l=}
\lfloor3.2d_{in}-1.4n_0\rfloor, & if  Case A, \nonumber\\
  \lfloor1.6n_0+ 3.2d_{in}\rfloor, & if  Case B.\label{dl-def}
\end{numcases}
By (\ref{lower-din-smaller}) and (\ref{n0-upper-bound}), we see that
\begin{numcases}{d_l\le}
3.2\alpha_2|V_2|, & if  Case A, \nonumber\\
 6.4\alpha_2|V_2|, & if  Case B.\label{dl-upper}
\end{numcases}

Then there is a vertex $v_{\scr 2}^{\scr *}$ in $V_2'\cup V_{12}^0$ of degree at least $|V_1'|-d_l.$
We will fix this vertex in what follows.
In fact, such a vertex $v_{\scr 2}^{\scr *}$ is in $V_2'$ by the facts that
\begin{equation}\label{v1dl}
    \delta(V_{12}^0,V_1')<(1-3\alpha_1)|V_1|\quad \mbox{and} \quad |V_1'|-d_l\geq (1-3\alpha_1)|V_1|,
\end{equation}
where $|V_1'|-d_l\geq (1-3\alpha_1)|V_1|$ holds because of (\ref{dl-upper}).

\noindent \textbf{Step 3. Finding a matching $M$ within $G[\Gamma(v_{\scr 2}^{\scr *},V_1')]$}

In this step, if $e_{in}\geq 1,$ we first find a matching within $G[V_1']$ of size at least $e_{in}/(2\bigtriangleup(G[V_1']))$.
We assume this by giving the following lemma.

\begin{LEM}
If $G$ is a graph with maximum degree $\Delta$, then $G$ contains a matching of size at least $\frac{|E(G)|}{2\Delta}$.
\end{LEM}

\pf We use induction on $|V(G)|$.
We may assume that the graph is connected. For otherwise, we are done by the induction hypothesis.
Let $e=xy\in E(G)$ be an edge and $G'=G-\{x,y\}$.
Since $|N_G(x)\cup N_G(y)|-|\{x,y\}|\le 2(\Delta -1)$, we have
$$e(G')\geq e(G)-2(\Delta-1)-1\geq e(G)-2\Delta.$$
Hence, by the induction hypothesis, $G'$ has a matching of size at least $\frac{e(G)-2\Delta}{2\Delta}=\frac{e(G)}{2\Delta}-1.$
Adding $e$ to the matching obtained in $G'$gives a matching of size at least $\frac{e(G)}{2\Delta}$ in $G$.
\qed

In case A, we take a matching in $G[V_1']$ of size at least $\max\{ \lfloor 11 d_{in}\rfloor, 11n_0\}$.
This is possible because
$$ \frac{e_{in}}{2\bigtriangleup(G[V_1'])}\geq \frac{e_{in}}{2\alpha_1|V_1'|/3}=\frac{3d_{in}}{2\alpha_1}\geq 11d_{in}$$
provided that $\alpha \leq (\frac{3}{22})^3$,
and
\begin{eqnarray}
      2e_{in} &\geq& |V_1'|\lceil(2n+3)/5 \rceil- |V_1'||V_2'|-(1-3\alpha_1)|V_1||V_{12}^0|\nonumber \\
       &\geq & |V_1'|\lceil(2n+3)/5 \rceil- |V_1'|(\lceil(2n+3)/5\rceil-n_0-|V_{12}^0|)-|V_1||V_{12}^0|+3\alpha_1|V_1||V_{12}^0|\nonumber \\
      &\geq & |V_1'|n_0+3\alpha_1|V_1||V_{12}^0|
\end{eqnarray}
implying that
$$
\frac{e_{in}}{2\bigtriangleup(G[V_1'])}\geq \frac{e_{in}}{2\alpha_1|V_1'|/3}\geq
\frac{|V_1'|n_0/2}{2\alpha_1|V_1'|/3}\geq\frac{3n_0}{4\alpha_1}\geq 11n_0
$$
provided that $\alpha \leq (\frac{3}{44})^3$.

By (\ref{dl-def}), $|V_1'|-\Gamma(v_{\scr 2}^{\scr *},V_1')\leq  d_l\le \lfloor 3.2 d_{in} \rfloor$,
we can then choose a matching $M$ from $\Gamma(v_{\scr 2}^{\scr *},V_1')$ such that
\begin{equation}\label{Msize}
 |M|=\max\{ \lfloor 7d_{in}\rfloor,  7 n_0 \}.
\end{equation}

In case B,
we take a matching in $G[V_1']$ of size at least $\lfloor 8d_{in} \rfloor$.
This is possible as
$$
\frac{e_{in}}{\bigtriangleup(G[V_1'])}\geq \frac{e_{in}}{2\alpha_1|V_1'|/3}=\frac{3d_{in}}{2\alpha_1}\geq \lfloor 8d_{in} \rfloor
$$
provided that $\alpha\le (\frac{3}{16})^3$.

By the second equality of (\ref{dl-def}),
$|V_1'|-\Gamma(v_{\scr 2}^{\scr *},V_1')\leq\lfloor  3.2d_{in}  +1.6n_0\rfloor$.
If $n_0<2d_{in}$, then $\lfloor 3.2d_{in}  +1.6n_0 \rfloor\leq \lfloor 7d_{in} \rfloor$.
Thus, there is a matching $M$ within $\Gamma(v_{\scr 2}^{\scr *},V_1')$ such that
\begin{numcases}{|M|=}
   \lfloor d_{in}, \rfloor & if  $n_0<2d_{in}$,  \nonumber\\
   0, & if $n_0\geq  2d_{in} $.  \label{M-size_B}
\end{numcases}
We fix $M$ for denoting the matching we defined in this step hereafter.

\noindent \textbf{Step 4. Insertion}

In this step, we insert vertices in $V_{12}^0$ into $V_1'- V(M)$.
Let $I=V_{12}^{0}=\{x_1,x_2,\cdots, x_{I}\}$, $U_1=\Gamma(v_{\scr 2}^{\scr *},V_1')-V(M)$,
and $U_2=V_2'$.
Then
(i)
\begin{eqnarray*}
      \delta(I, U_1) &\geq&  \delta (I, V_1')-|V(M)|-|V_1'-\Gamma(v_{\scr 2}^{\scr *},V_1')|\\
     &\geq & \alpha_1|V_2|/6-\max\{\lfloor 7d_{in} \rfloor, 7n_0\}-\lfloor 1.6n_0+3.2d_{in} \rfloor,\\
    &\geq & \alpha_1|V_2|/6-20.4\alpha_2|V_2|\quad (\mbox{by} \, (\ref{lower-din-smaller})\, \mbox{and}\,(\ref{n0-upper-bound}))\\
   &\geq & 3\alpha_2|V_2|\geq 3|I|\quad (\mbox{provided that}\, 23.4\alpha_2\le \alpha_1/6),
 \end{eqnarray*}
and from~(\ref{degree-V1'toV2'}), we have (ii)
\begin{equation*}
 \delta(U_1,U_2-\{v_{\scr 2}^{\scr *}\})\geq \lceil(2n+3)/5\rceil- \alpha_{1}|V_2|/3-1 > \alpha_2|V_2|\geq |I|.
\end{equation*}

By condition (i), there is a claw-matching $M_1$ between $I$ and $U_{1}$ centered in $I$.
Suppose that $\Gamma(x_{i}, M_1)= \{x_{i0}, x_{i1}, x_{i2}\}$.
We denote by $P_{x_i}$ the path $x_{i1}x_{i}x_{i2}$.
By (ii), there is a matching $M_2$ between $\{x_{i0}\,|\, 1\leq i\leq |I|\}$
 and $U_{2}-\{v_{\scr 2}^{\scr *}\}$ covering $\{x_{i0}\,|\, 1\leq i\leq |I|\}$.
 So far, we get two matchings $M_1$ and $M_2$.

Next we delete three types of edges not contained in
$$
(\bigcup_{i=1}^{|I|} E(P_{x_i}))\cup \{x_{i}x_{i0}\,:\,1\leq i\leq |I|\}.
$$
Those edges include edges incident to a vertex in $I$, edges incident to a vertex in
$$
\bigcup_{i=1}^{|I|}\left((\Gamma(x_{i1})-\Gamma(x_{i2}))\cup (\Gamma(x_{i2})-\Gamma(x_{i1}))\right),
$$
and one edge from the two edges connecting
a vertex in $\Gamma(x_{i1})\cap \Gamma(x_{i2})$ to both $x_{i1}$ and $ x_{i2}$, for each $i= 1, 2, \cdots, |I|$.

 For the resulting graph after the deletion of edges above, we contract each path $P_{x_i}\,(1\leq i\leq |I|)$
into a single vertex $v_{x_i}$. We call each $v_{x_i}$ a {\it wrapped vertex} and call $P_{x_i}$ the {\it preimage} of $v_{x_i}$.
Denote by $G^*$ the graph obtained by deleting and contracting the same edges as above, and let $U_2^*=V_2'$ and $U_1^*=V(G^*)-U_2^*$.
(We will need the following degree condition in the end of this proof.)
Since $|U^*_2|=|V_2'|\le (2/5+\alpha)n$, combining with~(\ref{degree-V1'toV2'}), we have
$$
deg(v_{x_i}, U^*_2)\geq |\Gamma(x_{i1}, U_2^*)\cap \Gamma(x_{i2}, U_2^*)|-1\ge 2n/5-\alpha_1|V_2|.
$$
By the above inequality and~(\ref{degree-V1'toV2'}), we get the first inequality below in~(\ref{G*-degree}).
Since one edge from the two edges connecting
a vertex in $\Gamma(x_{i1})\cap \Gamma(x_{i2})$ to both $x_{i1}$ and $ x_{i2}$ is deleted
 in $G^*$ for each $i= 1, 2, \cdots, |I|$, combining with~(\ref{v2V1 mindeg}), we have the second inequality as follows.
\begin{eqnarray}\label{G*-degree}
  \delta(U_{1}^*,U_{2}^*) &\ge & 2n/5 -\alpha_1|V_2|,\nonumber \\
  \delta(U_{2}^*,U_{1}^*) &\ge  & \delta(V_2',V_1')-1\ge (1-3\alpha_1)|V_1|-1.
\end{eqnarray}

Let $U'_{1}$ and $U'_{2}$ be the corresponding sets of $U_{1}$ and $U_{2}$, respectively, after
the contraction.
Let $T_W$ be the graph with
$$
V(T_{W})= \{x_{i0}, v_{x_i}\,:\, 1\leq i\leq |I|\} \cup (V(M_2)\cap U_{2}) \quad \mbox{and}
\quad
E(T_{W})= \{x_{i0}v_{x_i}\,:\, 1\leq i\leq |I|\} \cup E(M_2).
$$
By the construction,
$$
|V(T_{W})\cap U'_{1}|=|\{x_{i0}, v_{x_i}\,:\, 1\leq i\leq |I|\}|=2|I|,\,
|L(T_{W})\cap U'_{1}|=|\{v_{x_i}\,:\, 1\leq i\leq |I|\}|=|I|, \quad \mbox{and}
$$
$$
|V(T_{W})\cap U'_{2}|=|L(T_{I})\cap U'_{2}|=|V(M_2)\cap U'_{2}|=|I|.
$$

Notice that $T_{W}$ is a forest with $|I|$ components
and each vertex $x_{i0}\, (1\leq i\leq |I|)$ has degree 2 in $T_{W}$.
(We will make $T_{W}$ connected in the end by connecting each  $x_{i0}$ to $v_{\scr 2}^{\scr *}$.)
See a depiction of  this operation with $|I|=1$ in Figure~\ref{TI} below.

\begin{figure}[!htb]
\psfrag{x10}{$x_{10}$} \psfrag{x11}{$x_{11}$}
\psfrag{x12}{$x_{12}$} \psfrag{x1}{$x_{1}$} \psfrag{vx1}{$v_{x_1}$}
\begin{center}
  \includegraphics[scale=0.4]{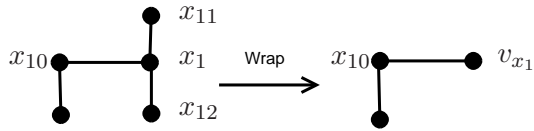}\\
\end{center}
\vspace{-3mm}
  \caption{$T_W$ with $|I|=1$}\label{TI}
\end{figure}


Let  $U_{I}^1=(V_1'-U_1)\cup U_1'-V(T_{ W})$, $U_{I}^2=U_2'-V(T_{ W})$,  and $G_I$ the resulting graph with $V(G_I)=U_{I}^1\cup U_{I}^2$.
We have that
\begin{eqnarray}\label{degree-UI}
    |U_{I}^1|&=&|V_1'|-3|I|=|V_1'|-3n_{12}^0, \quad  |U_{I}^2|=|V_2'|-|I|=|V_2'|-n_{12}^0,\nonumber\\
  \delta(U_{I}^1,U_{ I}^2)& \ge &\delta(V_1',V_2')-n_{12}^{ 0} \geq \lceil(2n+3)/5 \rceil -\alpha_1|V_2|/3-n_{12}^{ 0},\nonumber \\
  \delta(U_{I}^2,U_{I}^1)& \ge & \delta(V_2',V_1')-3n_{12}^{0}\geq (1-3\alpha_1)|V_1|-3n_{12}^{0}.
\end{eqnarray}

\noindent \textbf{Step 5. Matching Extension}

In this step, in the graph $G_I$,  we do some operation on the matching $M$ found in
Step 3.  Notice that the vertices in $M$ are unused in Step 4.
Recall that  $|M|\leq \max\{7n_0, \lfloor 7d_{in}\rfloor\}$.
By  $\lfloor d_{in}\rfloor\leq \alpha_2|V_2|$ from~(\ref{lower-din-smaller}) and
$n_0 \leq 2\alpha_2|V_2|$ from~(\ref{n0-upper-bound}), we get
\begin{equation}\label{EM-real}
   |M|\leq 14 \alpha_2|V_2|.
\end{equation}

\begin{figure}[!htb]
\psfrag{U1}{$U_1$} \psfrag{U2}{$U_2$}
\begin{center}
  \includegraphics[scale=0.4]{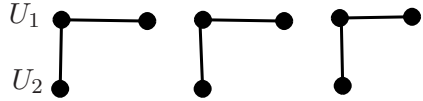}\\
\end{center}
\vspace{-3mm}
  \caption{$F_M$ with $|M|=3$}\label{FM}
\end{figure}

Hence, $\delta(U_{I}^{1},U_{I}^{2}-\{v_{\scr 2}^{\scr *}\})\geq \lceil(2n+3)/5\rceil -\alpha_1|V_2|/3-n_{12}^{0}-1 \geq |M|$.
Let $V_M$ be the set of vertices containing exactly one end of each edge in $M$.
Then there is a matching $M'$
between $V_{ M}$ and $U_{2}-\{v_{\scr 2}^{\scr *}\}$ covering $V_{M}$.
Let  $F_{M}$ be a forest with
$$
V(F_{M})=V(M)\cup (V(M')\cap U_{2}) \quad \mbox{and} \quad E(F_{M})=E(M)\cup E(M').
$$
Notice that
$$
|V(F_{M})\cap U_{1}|=2|M|, \quad |L(F_{M})\cap U_{1}|=|V(M)-V_{M}|=|M|,
$$
$$
 |V(F_{M})\cap U_{2}|=|L(F_{ M})\cap U_{2}|=|M|.
$$
Notice that $F_{M}$ has $|M|$ components, and all vertices in $V_{M}$ has degree 2.
(We will make $F_{M}$   a HIT later on by connecting each vertex in $V_{M}$ to the vertex $v_{\scr 2}^{\scr *} \in U_{2}$)
See Figure~\ref{FM} for a depiction of $F_{M}$ with $|M|=3$.

Let
$$
U_{\scr M}^1= U_{\scr I}^1-V(F_{\scr M}) \quad \mbox{and} \quad  U_{\scr M}^2= U_{\scr I}^2-V(F_{\scr M}).
$$

Notice that
\begin{eqnarray}\label{UM2}
    |U_{M}^1|&=& |U_{I}^1|-2|M|=|V_1'|-3 n_{12}^{0}-2|M|, \nonumber \\
    |U_{M}^2|&=& |U_{I}^2|-|M|=|V_2'|-n_{12}^{0}-|M|,
\end{eqnarray}
and
\begin{eqnarray}\label{UM1-2-degree}
     \delta(U_{M}^1,U_{M}^2) &=& \lceil(2n+3)/5\rceil - \alpha_1|V_2|/3-n_{12}^{ 0}-|M|,\nonumber\\
          \delta(U_{M}^2,U_{ M}^1)&\geq & (1-3\alpha_1)|V_1|-3n_{12}^{0}-2|M|.
\end{eqnarray}

\noindent \textbf{Step 6. Distribute Remaining vertices in $U_{\scr M}^1-\Gamma(v_{\scr 2}^{\scr *},V_1')$}

Let
$$
$$
 We may assume $n_l \geq 1$. For otherwise, we skip this step.
 By (\ref{dl-upper}), we have
\begin{numcases}{n_l\leq }\label{nl-upper}
3.2\alpha_2|V_2|, &  Case A, \nonumber\\
  6.4\alpha_2|V_2|, &  Case B.
\end{numcases}

By  $n_{12}^{0}\leq \alpha_2|V_2|$ from~(\ref{size of V2-v2'}) and $|M|\le 14\alpha_2|V_2|$ from~(\ref{EM-real}), we have (i)
\begin{eqnarray}\label{mindreUM1UM2}
  \delta(U_{\scr M}^1,U_{\scr M}^2) & \geq &  \lceil(2n+3)/5\rceil- \alpha_1|V_2|/3-n_{12}^0-|M|
   \geq \lceil(2n+3)/5\rceil- \alpha_1|V_2|/3- 15\alpha_2|V_2|\nonumber\\
&\ge & (1-3\alpha)|V_2|- \alpha_1|V_2|/3- 15\alpha_2|V_2|\quad(\mbox{as $\lceil(2n+3)/5\rceil\ge (1-3\alpha)(2/5+\alpha)n$})\nonumber\\
   &\geq & (1-3\alpha- \alpha_1/3-15\alpha_2)|V_2|)
  \geq  (1-\alpha_1)|V_2| \quad (\mbox{provided}\, 3\alpha+15\alpha_2 \le 2\alpha_1/3)\nonumber\\
  &\ge & (1-\alpha_1)|U_M^2|.
\end{eqnarray}
By (\ref{UM2}) and (\ref{nl-upper}), we have (ii)
\begin{eqnarray}
              |U_{\scr M}^2|- 10\alpha_1|V_2|- \lceil n_l/10\rceil -1 &\geq & |V_2'|- n_{12}^0-|M|- 16\alpha_1|V_2|-0.64\alpha_2|V_2| -2 \nonumber\\
                   &\geq & (1-\alpha_2-14\alpha_2-10\alpha_1-0.64\alpha_2-|V_2|/2)|V_2|\nonumber\\
                  &\geq &(1-11\alpha_1)|V_2| \quad (\mbox{provided}\, 15.64\alpha_2+ |V_2|/2\le \alpha_1)\nonumber\\
                   & > & 0 \quad (\mbox{provided}\, 11\alpha_1<1). \nonumber
\end{eqnarray}
Let $U_{R}=U_{\scr M}^1-\Gamma(v_{\scr 2}^{\scr *},V_1')$ and
denote $\left\lceil \frac{|U_{R}|}{10}\right\rceil=l$.
Suppose first that $|U_{R}|\geq 2$.
We partition $U_{ R}=U_{R_1}\cup U_{R_2}\cup\cdots\cup U_{R_l}$
arbitrary such that each set contains at least 2 and at most $|U_R|/10$ vertices.
Then by the conditions (i) and (ii), for each $1\leq i\leq l$,
there is a vertex $y_{i} \in U_{2}-\{v_{\scr 2}^{\scr *}\}$ which is common to all vertices in $U_{R_i}$,
and is not used by any other $U_{R_j}$ with $j\neq i$.
Let $T_{R}$ be the graph with
$$
V(T_{R})=U_{R}\cup \{y_{i}\,:\,1\leq i\leq l \} \quad \mbox{and}\quad E(T_{R})= \{xy_{i}\,:\, x\in U_{R_i}, 1\leq i\leq l\}.
$$
Suppose now $|U_{\scr R}|=1$, let $U_{R}=\{x_{R}\}$. Choose $x'_{R} \in U_{M}^1-U_{R}$
and $y_{R} \in U_{M}^2-\{v_{\scr 2}^{\scr *}\}$ be a vertex common to $x_{R}$ and $x'_{R}$.
Let $T_{R}$ be a tree with
$$
V(T_{R})=\{x_{R}, x'_{R}, y_{R}\} \quad \mbox{and} \quad E(T_{R})=\{x_{R}y_{R}, x'_{R}y_{R}\}.
$$
By the construction,
$$
|V(T_{R})\cap  U_{M}^1|=|L(T_{R}\cap  U_{M}^1|=\max\{|U_{R}|,2\},\quad
|V(T_{ R})\cap U_{M}^2|=l,  \quad \mbox{and} \quad |L(T_{ R}\cap U_{M}^2|=0.
$$
Notice that $T_{R}$ is not connected when $|U_{R}|\geq 17$ and that $T_{R}$ may have degree 2 vertices in $V(T_R)\cap U_{M}^2$.
 Later on, by  joining each vertex in $T_{R}\cap U_{M}^2$  to a vertex of a tree,
we will make the resulting graph connected,  and thereby eliminating the possible degree 2 vertices in $T_{R}$.
Let
$$
U_{R}^1= U_{M}^1-V(T_{ R})\quad \mbox{and} \quad U_{R}^2= U_{M}^2-V(T_{R}).
$$
Then we have
\begin{eqnarray}\label{UR12-size}
               |U_{R}^1| &=& |U_{ M}^1|-n_l =|V_1'|-3n_{12}^0-2|M|-\max\{2, n_l\}, \nonumber\\
               |U_{R}^2|&=&|U_{M}^2|-\lceil n_l/10\rceil = |V_2'|-n_{12}^0- |M|-\lceil n_l/10\rceil,
\end{eqnarray}
and
\begin{eqnarray}\label{deltaUR21}
 \delta(U_{\scr R}^1,U_{\scr R}^2) &\geq &\lceil(2n+3)/5\rceil - \alpha_1|V_2|/3-n_{12}^0-|M|-\lceil n_l/10\rceil,\nonumber\\
 \delta(U_{\scr R}^2,U_{\scr R}^1)&=& (1-3\alpha_1)|V_1|-3n_{12}^0-2|M|-max\{2, n_l\}.
\end{eqnarray}
Let $G_R$ be the subgraph of $G$ induced on $U_{\scr R}^1\cup U_{\scr R}^2$.

\noindent \textbf{Step 7. Completion of a HIST in $G_R$}

 In this step, we find a HIST $T_{\scr main}$ in $G_R$ such that
\begin{eqnarray*}
  |L(T_{\scr main})\cap U_{ R}^1|=|L(T_{W})|/2+|L(F_{M})\cap U_{ I}^1|+ |L(T_{ R})\cap U^1_{M}| &=&  \\
  |L(T_{\scr main}\cap U_R^2)|=|L(T_{ W})|/2+|L(F_{ M})\cap U_{I}^2|+|L(T_{ R})\cap U^1_{ M}|. &&
\end{eqnarray*}
By the construction of $F_{ M}$ and $T_{\scr R}$, we have
$|L(F_{M})\cap U_{ I}^1|=|L(F_{ M})\cap U_{I}^2|$
and $|L(T_{R})\cap U_{M}^1|-|L(T_{R})\cap U_{M}^2|= \max\{2, n_l\}$,
respectively.
So
\begin{equation}\label{leaves-difference}
   |L(T_{\scr main})\cap U_{R}^2|-|L(T_{\scr main})\cap U_{R}^1|= \max\{2, n_l\}.
\end{equation}

 Notice that $v_{\scr 2}^{\scr *} \in  U_{ R}^2$, $v_{\scr 2}^{\scr *}$ is adjacent to each vertex in $U_R^1$,
 and $V_1'-\Gamma(v_{\scr 2}^{\scr *},V_1')\subseteq V(U^1_{\scr R})$.
We now construct $T_{\scr main}$ step by step.

\noindent  \textbf{Step 7.1 }

Let $T^1_{\scr main}$ be the graph with
$$
V(T^1_{\scr main})=\{v_{\scr 2}^{\scr *}\}\cup U_{\scr R}^1  \quad \mbox{and} \quad E(T^1_{\scr main})= \{v_{\scr 2}^{\scr *}x\,|\, x\in U_{\scr R}^1 \}.
$$
To make the requirement of (\ref{leaves-difference}) possible, we need to make at least
\begin{eqnarray}\label{d3}
  d_{3} &=&|U^1_{R}|-|U_{R}^2|+\max\{2, n_l\},\nonumber\\
   &=& |V_1'|-|V_2'|-2n_{12}^0-|M|+\lceil n_l/10\rceil
\end{eqnarray}
vertices in $U_{R}^1$ with degree at least 3 in $T_{\scr main}$, where the last inequality above follows from (\ref{UR12-size}).
Hereinafter,  we assume that $\max\{2, n_l\}=n_l$. Since the proof for $\max\{2, n_l\}=2$ follows the same idea, we skip the
details.

Since all vertices in $U_{R}^1 $ are included in $T^1_{\scr main}$ and $T^1_{\scr main}$ is connected,
each vertex in $T^1_{\scr main}$ needs to join to at least two distinct vertices from $U_{R}^2- \{v_{\scr 2}^{\scr *}\}$
to have degree no less than 3. Hence, to make  a desired HIST $T_{\scr main}$, it is necessary that
 \begin{eqnarray}\label{df*}
   d_{f*} &=& |U_{\scr R}^2|-1-2d_{3} \nonumber\\
    &=& |V_2'|-n_{12}^0-e_M-\lceil n_l/10\rceil-1-2d_{3} \nonumber\\
      &=& 3|V_2'|-2|V_1'|+3n_{12}^0+|M|-3 \lceil n_l/10\rceil-1\nonumber\\
      &\geq&0.
 \end{eqnarray}
We show (\ref{df*}) is true, separately,  for each of Case A and Case B. For Case A, notice that
$$
|V_1'|=n-\lceil(2n+3)/5\rceil+ n_0 \quad \mbox{and} \quad |V_2'|=\lceil(2n+3)/5\rceil-n_0-n_{12}^0.
$$
Hence,
$$
3|V_2'|=3\lceil(2n+3)/5\rceil-3n_0-3n_{12}^0 \quad \mbox{and} \quad  2|V_1'|=2n-2\lceil(2n+3)/5\rceil+2n_0.
$$
Thus,
\begin{eqnarray*}
   d_{f*}&=&5\lceil(2n+3)/5\rceil-2n-5n_0-3n_{12}^0+3n_{12}^0+|M|-3\lceil n_l/10\rceil-1 \\
   &\geq&2-5n_0+|M|-3\lceil\lfloor3.2d_{in}\rfloor/10\rceil\quad (\mbox{by}~n_l\le d_l \lfloor 3.2d_{in}-1.4n_0\rfloor\,\mbox{from~(\ref{dl-def})})\\
   &=&2-5n_0+\max\{7n_0, \lfloor7d_{in}\rfloor\}-3\lceil\lfloor3.2d_{in}\rfloor/10\rceil\\
   &\geq& 0.
\end{eqnarray*}
Now we show (\ref{df*}) is true for case B. Notice that
$$
|V_1'|=n-\lceil(2n+3)/5\rceil+ n_0 \quad \mbox{and} \quad |V_2'|=\lceil(2n+3)/5\rceil+ n_0-n_{12}^0.
$$
So
$$
3|V_2'|=3\lceil(2n+3)/5\rceil+3n_0-3n_{12}^0 \quad \mbox{and} \quad  2|V_1'|=2n-2\lceil(2n+3)/5\rceil+2n_0.
$$
Recall that $n_0\ge 1 $ in this case. We have
 \begin{eqnarray*}
    d_{f*}&=&5\lceil(2n+3)/5\rceil -2n + n_0-3n_{12}^0+3n_{12}^0+|M|-3\lceil n_l/10\rceil-1\nonumber \\
    &\geq & 2+ n_0+|M|-3\lceil n_l/10\rceil\nonumber \\
    &=& 2+ n_0+|M|-3\lceil \lfloor3.2d_{in}+1.6n_0\rfloor/10\rceil \quad (\mbox{by}~n_l\le \lfloor3.2d_{in}+1.6n_0\rfloor\,\mbox{from~(\ref{dl-def})})\nonumber\\
    &\ge & \left\{
             \begin{array}{ll}
               2+n_0+\lfloor d_{in}\rfloor- \lfloor9.2d_{in}/10\rfloor-\lfloor4.8n_0/10\rfloor-1\ge 0, & \hbox{if $n_0<2d_{in}$;} \nonumber\\
               2+n_0- \lfloor9.2d_{in}/10\rfloor-\lfloor4.8n_0/10\rfloor-1\ge 0, & \hbox{if $n_0\ge 2d_{in}$.}
             \end{array}
           \right.
 \end{eqnarray*}

We now in Step 2 below show that there is a way to make exactly $d_{f*}$ vertices in $T^1_{\scr main}$
with degree 3 by joining each to two distinct vertices from $U_{\scr R}^2-\{v_{\scr 2}^{\scr *}\}$.

\noindent \textbf{Step 7.2}

 We first take  $2d_{3}$ vertices from $U_{R}^2-\{v_{\scr 2}^{\scr *}\}$.
 For those $2d_{3}$ vertices, pair them up into $d_{3}$  pairs.
 We show that for each pair of vertices, they have at least $d_{3}$ common neighbors in $U_{R}^1.$
Using (\ref{deltaUR21}), $|M|\le14\alpha_2|V_2|$ from~(\ref{EM-real}), $n_l\le d_l\le 6.4\alpha_2|V_2|$ from~(\ref{dl-upper}),
we have
\begin{eqnarray}\label{deltaUR21-accurate}
  \delta(U_{ R}^2,U_{ R}^1)  &\ge & (1-3\alpha_1)|V_1|-3n_{12}^0-2|M|- \max\{2, n_l\}\nonumber \\
   &\ge & |V_1|-3\alpha_1|V_1|- 3\alpha_2|V_2|-28\alpha_2|V_2|-6.4\alpha_2|V_2|\nonumber\\
   &\ge & |V'_1|-|V_1-V_1|-37.4\alpha_2|V_2|- 3\alpha_1|V_1|.
\end{eqnarray}
Since $|U_R^1|\le |V'_1|$,  we know that any two vertices in $U_{R}^2$ have at least
 \begin{eqnarray}
   n_c &=& |V_1|-2|V_1'-V_1|-74.8\alpha_2|V_2|- 6\alpha_1|V_1|\nonumber \\
    &\geq & (3/5-\alpha)n-76.8\alpha_2|V_2|- 6\alpha_1|V_1|\quad (\mbox{by}\, |V_1'-V_1|=|V_0'|\le  |V_2-V_2'|\le \alpha_2|V_2|) \nonumber \\
    &\geq & 3n/5-10\alpha_1|V_1| \quad (\mbox{provided that}\, 76.8\alpha_2+3\alpha\le  4\alpha_1)\nonumber
 \end{eqnarray}
common neighbors in $U_R^1$.
On the other hand,
 \begin{eqnarray}\label{d3-upper}
  d_{3}&=& |V_1'|-|V_2'|-2n_{12}^0-|M|-\lceil n_l/10\rceil \nonumber \\
   &\leq &(3/5-\alpha)n-(2n/5-2\beta n-|V_2-V_2'|)+(1.6n_0+3.2\lfloor d_{in}\rfloor )/10+1\nonumber \\
   &=&n/5-\alpha n+2\beta n+|V_2-V_2'|+(3.2\alpha_2|V_2|+3.2\alpha_2|V_2|)/10\quad  (\mbox{by}~(\ref{lower-din-smaller}) \, \mbox{and}\, (\ref{n0-upper-bound}))\nonumber \\
  &\leq & n/5-\alpha n+2\beta n+\alpha_2|V_2|+0.64\alpha_2|V_2|\nonumber \\
  &\leq &n/5+2\alpha_1|V_2|<n_c\,(\mbox{provided}\, 12\alpha_1<2/5). \nonumber
 \end{eqnarray}
Denote by $\{u_1^1, u_1^2\},\{u_2^1, u_2^2\},\cdots,\{u_{d_{3}}^1, u_{d_{3}}^2\}$
the $d_{3}$ pairs of vertices from $U_{R}^2-\{v_{\scr 2}^{\scr *}\}.$
Then by the above argument, we can choose $d_{3}$ distinct vertices say $v_1,v_2, \cdots, v_{d_{3}}$
from $L(T^1_{\scr main})$
 such that $v_i\sim  u_i^1,  u_i^2$ for all $1\leq i\leq d_{3}$.

 Let $T^2_{\scr main}$  be the graph with
$$
V(T^2_{\scr main})=V(T^1_{\scr main})\cup \{u_i^1,u_i^2\,:\,1\leq i\leq d_{3} \}\quad \mbox{and}\quad
E(T^2_{\scr main})=E(T^1_{\scr main})\cup \{v_iu_i^1,v_iu_i^2\,:\,1\leq i\leq d_{3} \}.
$$
If $V(G_R)-V(T^2_{\scr main})=\emptyset$, we let $T_{\scr main}=T^2_{\scr main}$.
For otherwise, we need one more step to finish constructing $T_{\scr main}$.

\noindent \textbf{Step 7.3}

For the remaining vertices in $U_{\scr R}^2-V(T^2_{\scr main}) $,
we show that each of them has a neighbor
in $S(T^2_{\scr main})\cap U_{\scr R}^1$; that is,  a neighbor in $U_{\scr R}^1$
of degree 3 in $V(T^2_{\scr main})$. This is clear, as by (\ref{deltaUR21-accurate}),  we have
\begin{eqnarray*}
  \delta(U_{\scr R}^2,U_{\scr R}^1)& \geq & |V_1'|-|V_1'-V_1|-37.4\alpha_2|V_2|- 3\alpha_1|V_1|\\
  &\geq&  |U_R^1|-38.4\alpha_2|V_2|- 3\alpha_1|V_1|\quad (\mbox{by}\, |V_1'-V_1|\le |V_2-V_2'|\le \alpha_2|V_2|).
\end{eqnarray*}
Since $|S(T^2_{\scr main})\cap U_{\scr R}^1|=d_3$, and
\begin{eqnarray*}
  d_{3} &=& |V_1|-|V_2'|-2n_{12}^0-2|M|+\lceil n_l/10\rceil \\
  &\geq&(3/5-\alpha)n-(2/5+\alpha)n-2\alpha_2|V_2|-28\alpha_2|V_2|+0.64\alpha_2|V_2| \\
    &\geq& n/5-2\alpha n-29.36\alpha_2|V_2|\\
    &>&38.4\alpha_2|V_2|+ 3\alpha_1|V_1| \quad (\mbox{provided}\, 2\alpha+67.76\alpha_2+3\alpha_1<1/5).
\end{eqnarray*}
Now, we join an edge between each vertex in $U_{R}^2-V(T^2_{\scr main})$ and a neighbor of the
 vertex in $S(T^2_{\scr main})\cap U_{R}^1$.
Let $T_{\scr main}$  be the resulting tree.  By the construction procedure,
it is easy to verify that $T_{\scr main}$
is a  HIST of $G_{R}$.

\noindent \textbf{Step 8. Connecting $T_{W}$, $F_{M}$, $T_{R}$,  and $V(T_{\scr main})$ into a connected graph}

In this step, we connect  $T_{ W}$, $F_{M}$, $T_{R}$, and $V(T_{\scr main})$ into  a connected graph.
Recall that each degree 2 vertex in $T_{W}$ and  $F_{ M}$ is a neighbor of $v_{\scr 2}^{\scr *}$.
We join an edge connecting $v_{\scr 2}^{\scr *}$ in $V(T_{\scr main})$ and each degree 2 vertex in $T_{W}$ and  $F_{M}$.
By the argument  in step 7.3 above, we know each vertex in $V(T_{R})\cap U_{M}^2$
has a neighbor in $S(T_{\scr main})\cap U_{R}^1$.
Thus, we join an edge between each vertex in $V(T_{R})\cap U_{M}^2$ to exactly one of its neighbor
in $S(T_{\scr main})\cap U_{R}^1$.
Let $T^*$ be the final resulting graph.
Notice that  $I=V_{12}^{0}=\{x_1,x_2,\cdots, x_{I}\}\subseteq L(T^*)$ is the
set of the wrapped vertices from Step 4.
Recall that $G^*$ is the  graph obtained from $G$ be deleting and contracting edges from Step 4.
Then by the constructions of $T_{W}$, $F_{M}$, $T_{R}$,
and $T_{\scr main}$,
we see that $T^*$ is a HIST of $G^*$ with $|L(T^*)\cap U_1^*|=|L(T^*)\cap U_2^*|$.

\noindent \textbf{Step 9. Finding a cycle on $L(T^*)$}

Denote
$$
 U_L^1=L(T^*)\cap U_1^*,\,\, U_L^2=L(T^*)\cap U_2^* \quad \mbox{and}\quad  G_L=G[E_G(U_L^1, U_L^2)].
$$
Notice that $G_L$ is a balanced bipartite graph. And
\begin{eqnarray*}
  |S(T^*)\cap U_1^*| &=& d_3\le n/5+2\alpha_1|V_2|\quad (\mbox{by}\, (\ref{d3-upper})) \\
   |S(T^*)\cap U_2^*| &=& 1+\lceil n_l/10\rceil \le 2+0.64\alpha_2|V_2| \quad (\mbox{by $n_l\le d_l\le 6.4\alpha_2|V_2|$ from~(\ref{dl-upper})}).
\end{eqnarray*}
Thus by (\ref{G*-degree}),
\begin{eqnarray*}
  \delta_{G^*}(U_{L}^1,U_{ L}^2) &\ge &2n/5 -\alpha_1|V_2|-(2+0.64\alpha_2|V_2|)>3n/10>|U_L^2|/2+1, \\
  \delta_{G^*}(U_{L}^2,U_{L}^1) &\ge  &  (1-3\alpha_1)|V_1|-1-(n/5+2\alpha_1|V_2|)>n/3>|U_L^1|/2+1.
\end{eqnarray*}
By Lemma~\ref{bi-Hamiltonian},
$G_L$ contains a hamiltonian cycle $C'$.

\noindent \textbf{Step 10. Unwrap vertices in $V(C')\cap \{v_{x_1},v_{x_2},\cdots, v_{x_{|I|}}\} $}

On $C'$, replace each vertex $v_{x_i}$ with its preimage $P_{x_{i}}=x_{i1}x_{i}x_{i2}$ for each $i=1, 2, \cdots, |I|$.
Denote the resulting cycle by $C$. Recall that $x_{i1}, x_{i2}\in \Gamma(v_{\scr 2}^{\scr *})$ by the
choice of $x_{i1}$ and $ x_{i2}$.
Let $T$ be the graph on $V(G)$ with
$$
E(T)=E(T^*)\cup\{v_{\scr 2}^{\scr *}x_{i1}, v_{\scr 2}^{\scr *}x_{i2}\,:\,i=1, 2, \cdots, |I|\}.
$$
Then $T$ is a HIST of $G$.
Let $H=T\cup C$. Then $H$ is an SGHG of $G$.

The proof of Extremal Case 2 is finished.
\qqed

\bibliographystyle{plain}
\bibliography{SSL-BIB}

\end{document}